\documentclass[10pt]{amsart2000}
\usepackage{amsmath2000,amssymb,amscd2000,latexsym}

\newcommand{\C}{{\mathbb{C}}}

\newcommand{\N}{{\mathbb{N}}}

\newcommand{\R}{{\mathbb{R}}}
\newcommand{\T}{{\mathbb{T}}}
\newcommand{\Z}{{\mathbb{Z}}}

\newcommand{\Mf}{{\mathfrak M}}

\newcommand{\Ah}{{\mathcal A}}
\newcommand{\Bh}{{\mathcal B}}
\newcommand{\Ch}{{\mathcal C}}
\newcommand{\Dh}{{\mathcal D}}
\newcommand{\Eh}{{\mathcal E}}

\newcommand{\Jh}{{\mathcal J}}

\newcommand{\Ph}{{\mathcal P}}

\newcommand{\Uh}{{\mathcal U}}

\newcommand{\Zh}{{\mathcal Z}}

\newcommand{\be}{\mathbf{1}}
\newcommand{\bi}{\bar{\imath}}

\newcommand{\cpr}{\mathrm{cpr }\,}
\newcommand{\card}{\mathrm{card}}
\newcommand{\diam}{\mathrm{diam}}
\newcommand{\dist}{\mathrm{dist}}
\newcommand{\dr}{\mathrm{dr}\,}
\newcommand{\ev}{\mathrm{ev }}

\newcommand{\halb}{\frac{1}{2}}

\newcommand{\id}{\mathrm{id}}
\newcommand{\ord}{\mathrm{ord}\,}

\newcommand{\Prim}{\mathrm{Prim}}

\newcommand{\Sd}{\mathrm{Sd}}

\newcommand{\verk}{\mbox{\scriptsize $\,\circ\,$}}

\newcounter{number}[subsection]
\newcounter{altnumber}[section]

\newenvironment{nummer}{\refstepcounter{number}{\noindent\arabic{section}.\arabic{subsection}.\arabic{number}}}{}
\newenvironment{altnummer}{\refstepcounter{altnumber}{\noindent\arabic{section}.\arabic{altnumber}}}{}

\newcommand{\bn}{\noindent\begin{nummer} \rm}
\newcommand{\en}{\end{nummer}}

\newcommand{\altbn}{\noindent \begin{altnummer} \rm}
\newcommand{\alten}{\end{altnummer}}

\newenvironment{ntheorem}{\noindent {\sc Theorem:} \it}{}
\newenvironment{nlemma}{\noindent {\sc Lemma:} \it}{}
\newenvironment{nprop}{\noindent {\sc Proposition:} \it}{}
\newenvironment{ndefn}{\noindent {\sc Definition:} \it}{}
\newenvironment{ncor}{\noindent {\sc Corollary:} \it}{}

\newenvironment{nremark}{\noindent {\sc Remark:} }{}
\newenvironment{nremarks}{\noindent {\sc Remarks: }}{}

\newenvironment{nproof}{\noindent {\sc Proof:}}{\mbox{}\hfill \rule[-.2ex]{.25em}{1.8ex}}

\begin{document}

\title{{\sc Decomposition Rank of Subhomogeneous $C^*$-Algebras}}

\author{Wilhelm Winter}
\address{Mathematisches Institut der Universit\"at M\"unster\\ 
Einsteinstr. 62\\ D-48149 M\"unster}
\curraddr{{\sc Department of Mathematics\\ Texas A{\&}M University\\
 College Station\\
Texas\\
USA}}

\email{wwinter@math.uni-muenster.de}

\date{October 2002}
\thanks{\it Partially supported by the Deutsche Forschungsgemeinschaft}

\setcounter{section}{-1}


\begin{abstract}
We analyze the decomposition rank (a notion of covering dimension for nuclear $C^*$-algebras introduced by E.\ Kirchberg and the author) of subhomogeneous $C^*$-algebras. In particular we show that a subhomogeneous $C^*$-algebra has decomposition rank $n$ if and only if it is recursive subhomogeneous of topological dimension $n$ and that $n$ is determined by the primitive ideal space. \\
As an application, we use recent results of Q.\ Lin and N.\ C.\ Phillips to show the following: 
Let $A$ be the crossed product $C^*$-algebra coming from a compact smooth manifold and a minimal diffeomorphism. Then the decomposition rank of $A$ is dominated by the covering dimension of the underlying manifold.
\end{abstract}

\maketitle

\section{Introduction}

In \cite{KW}, E.\ Kirchberg and the author introduced the decomposition rank; this is a noncommutative generalization of topological covering dimension. If $A$ is a nuclear $C^*$-algebra, the decomposition rank of $A$, $\dr A$, is defined by imposing a certain condition on systems of completely positive (c.p.) approximations of $A$; see Section 1 for details.\\
It may happen that $A$ has some obvious underlying topological space $X$; in this case it is natural to ask wether $\dr A$ is related to the covering dimension of $X$, $\dim X$, in any way. There are several candidates of such underlying spaces that come to mind, like the spectrum $\hat{A}$ or the primitive ideal space $\Prim \, A$, but there might also be some space $X$ involved in the construction of $A$, e.g.\ if $A$ is the $C^*$-algebra generated by a group action on $X$.\\
The decomposition rank behaves very well if $A$ is a continuous trace algebra, for in this case $\hat{A}$ is a locally compact Hausdorff space and we have $\dr A = \dim \hat{A}$ by \cite{KW}, Proposition 3.12. For more general type $I$ $C^*$-algebras the situation is less obvious. In these notes we are mainly concerned with the case where $A$ is subhomogeneous, i.e.\ has irreducible representations of dimension at most $N$ for some $N \in \N$. Among these algebras, recursive subhomogeneous algebras of finite topological dimension (introduced in \cite{Ph}) are particularly tractable. Such an algebra can be written as an iterated pullback of algebras of the form $\Ch(X_j) \otimes M_{r_j}$ with finite-dimensional compact spaces $X_j$. In this case, the topological dimension of $A$ coincides with $\max_k \{\dim (\Prim_k A)\}$, where $\Prim_k A$ is the locally compact space of kernels of $k$-dimensional irreducible representations, as follows from a theorem of Phillips. Now if $A$ is subhomogeneous with finite decomposition rank, then the same theorem implies that $A$ in fact is recursive subhomogeneous of topological dimension at most $\dr A$. The aim of the present article is to prove a converse, namely that, if $A$ is recursive subhomogeneous of topological dimension $n$, then $\dr A \le n$.\\
Recursive subhomogeneous algebras play an important r\^ole in Elliott's classification program. Roughly speaking, the Elliott conjecture says that separable, simple, stably finite, nuclear $C^*$-algebras are classified by their $K$-theory data (cf.\ \cite{Ro}, Conjecture 2.2.5). If this was true, it would follow from theorems of Elliott and Thomsen about the range of the invariant that any such $C^*$-algebra (provided $K_0$ is weakly unperforated) is an inductive limit of recursive subhomogeneous algebras of toplogical dimension at most 2 (cf.\ \cite{Ro}, Theorem 3.4.4, and Example \ref{range} below). \\
By recent work of Lin and Phillips certain crossed product $C^*$-algebras also admit such a direct limit decomposition. More precisely, let $M$ be a compact smooth manifold, $\alpha : M \to M$ a minimal diffeomorphism and $A := \Ch(M) \times_\alpha \Z$ the crossed product $C^*$-algebra. Then $A$ can be written as an inductive limit of recursive subhomogeneous algebras of topological dimension at most $\dim M$. As a consequence we see that $\dr A \le \dim M$. A special case of this phenomenon already occurred in a theorem of Elliott and Evans which says that irrational rotation algebras are limit circle algebras (and thus have decomposition rank one). At the present stage this setting seems to be the only systematic way to obtain information on the decomposition rank of crossed products, since in general it is very hard to construct c.p.\ approximations for crossed products with sufficiently good properties.

The paper is organized as follows: In Section 1 we recall the definition of the decomposition rank and of recursive subhomogeneous algebras.  Furthermore, we state our main result (Theorem \ref{rec.subh.}), namely that a unital subhomogeneous algebra has finite decomposition rank $n$ iff it is recursive subhomogeneous of topological dimension $n$, and prove its easy part; we also give a nonunital version. Furthermore, we consider various examples, including the noncommutative $CW$-complexes of \cite{Ped} and $C^*$-algebras of minimal diffeomorphisms.\\
The remaining sections are devoted to the proof of the difficult part of the theorem. Since our argument is quite complicated, we first describe the ideas in the commutative setting (Section 2); the proof of the general case will be modelled after this outline. In Section 3 we deduce a lifting result for centers of certain subhomogeneous $C^*$-algebras. This is used in Section 4 to obtain an approximate lifting result for so-called piecewise commuting maps. In Section 5 we develop a topological concept which might be called relative barycentric subdivision. The actual proof of Theorem \ref{rec.subh.} is given in Section 6.    
  
We would like to thank J.\ Cuntz, S.\ Echterhoff, E.\ Kirchberg, N.\ C.\ Phillips and W.\ Werner for many helpful comments and discussions.

\section{Decomposition rank and recursive subhomogeneous $C^*$-algebras}

\altbn{\label{d-dr}}
Recall from \cite{Tak} that nuclear $C^*$-algebras are characterized by the completely positve approximation property, i.e.\ $A$ is nuclear if and only if there is a net of finite-dimensional $C^*$-algebras $F_\lambda$ and completely positive contractive (c.p.c.) maps $A \stackrel{\psi_\lambda}{\longrightarrow} F_\lambda \stackrel{\varphi_\lambda}{\longrightarrow} A$ such that $\varphi_\lambda \verk \psi_\lambda$ converges to $\id_A$ pointwise. We then say $(F_\lambda,\psi_\lambda,\varphi_\lambda)_\Lambda$ is a system of c.p.\ approximations for $A$. Based on this approximation property, one may define a noncommutative version of covering dimension as follows:

\begin{ndefn}
(cf.\ \cite{KW}, Definitions 2.2 and 3.1) Let $A$ be a separable $C^*$-algebra. \\
(i) A c.p.\ map $\varphi : \bigoplus_{i=1}^s M_{r_i} \to A$ has strict order zero, $\ord \varphi = 0$, if it preserves orthogonality, i.e., $\varphi(e) \varphi(f) = \varphi(f) \varphi(e) = 0$ for all $e,f \in \bigoplus_{i=1}^s M_{r_i}$ with $ef = fe = 0$.\\ 
(ii) A c.p.\ map $\varphi : \bigoplus_{i=1}^s M_{r_i} \to A$ is $n$-decomposable, if there is a decomposition $\{1, \ldots, s\} = \coprod_{j=0}^n I_j$ s.t.\ the restriction of $\varphi$ to $\bigoplus_{i \in I_j} M_{r_i}$ has strict order zero for each $j \in \{0, \ldots, n\}$; we say $\varphi$ is $n$-decomposable w.r.t.\ $\coprod_{j=0}^n I_j$.\\
(iii) $A$ has decomposition rank $n$, $\dr A = n$, if $n$ is the least integer such that the following holds: Given $\{b_1, \ldots, b_m\} \subset A$ and $\varepsilon > 0$, there is a c.p.\ approximation $(F, \psi, \varphi)$ for $b_1, \ldots, b_m$ within $\varepsilon$ (i.e., $\psi:A \to F$ and $\varphi:F \to A$ are c.p.c.\ and $\|\varphi \psi (b_i) - b_i\| < \varepsilon$) such that $\varphi$ is $n$-decomposable. If no such $n$ exists, we write $\dr A = \infty$.  
\end{ndefn}

This notion is a variation of the completely positive rank ($\cpr A$), which was introduced in \cite{Wi1}. It has good permanence properties; for example, it behaves well with respect to quotients, inductive limits, hereditary subalgebras, unitization and stabilization.\\
Both ranks generalize topological covering dimension, i.e., if $X$ is a locally compact second countable space, then $\cpr \Ch_0(X) = \dr \Ch_0(X) = \dim X$; see \cite{KW} for details. 
\alten

\altbn{\label{order-zero}}
In \cite{Wi1}, Proposition 4.4.1(a), maps of strict order zero were characterized as follows: \\
If $\varphi : F \to A$ is c.p.c.\ with $\ord \varphi = 0$, then there is a unique $*$-homomorphism $\pi_\varphi : CF \to A$ such that $\pi_\varphi(g \otimes x) = \varphi(x) \; \forall \, x \in F$, where $CF$ is the cone $\Ch_0((0,1]) \otimes F$ over $F$ and $g := \id_{(0,1]}$ is the canonical generator of $\Ch_0((0,1])$. Conversely, any $*$-homomorphism $\pi : CF \to A$ induces such a c.p.c.\ map $\varphi$ of strict order zero.\\
The $*$-homomorphism $\pi_\varphi$ extends to a $*$-homomorphism $\pi_\varphi'' : (CF)'' \to A''$ of von Neumann algebras. Then we have $\varphi (x) = \varphi(\be_F) \sigma (x) = \sigma(x) \varphi(\be_F) \; \forall \, x \in F$, where $\sigma : F \to A''$ is the $*$-homomorphism coming from the composition of the natural unital embedding $F \hookrightarrow (CF)''$ and $\pi_\varphi''$. Note that each $h \in C^*(\varphi(\be_F))$, $0 \le h \le \be$, defines a c.p.c.\ map $\hat{\varphi}: F \to A$ by $\hat{\varphi}(\,.\,):= h \sigma(\, .\,)$ s.t.\ $\ord \hat{\varphi} = 0$ and $\| \hat{\varphi} - \varphi\| = \|h - \varphi(\be_F)\|$. 
\alten

\altbn{\label{d-rec.subh.}}
\begin{ndefn}
(cf.\ \cite{Ph}, Definition 1.1) A recursive subhomogeneous algebra is a unital $C^*$-algebra defined recursively as follows:\\
(1) $M_r$ is a recursive subhomogeneous algebra for any $r \in \N$.\\
(2) If $B$ is a recursive subhomogeneous algebra, $\Omega$ is a compact Hausdorff space, $X \subset \Omega$ a closed subspace, $r \in \N$ and $\pi : B \to \Ch(X) \otimes M_r$ a unital $*$-homomorphism, then $A:= B \oplus_{\pi,X} (\Ch(\Omega) \otimes M_r) := \{ (b,f) \in B \oplus (\Ch(\Omega) \otimes M_r) \, | \, \pi(b)(t) = f(t) \; \forall \, t \in X \}$ is a recursive subhomogeneous algebra.
\end{ndefn}

We see from the definition that any algebra $A$ as above can be written as an iterated pullback involving base spaces $\Omega_k$; the topological dimension of such a decomposition is then defined as $\max_k \{\dim \Omega_k\}$. However, the decomposition is highly nonunique as easy examples show; by definition, the topological dimension of the algebra $A$ is the least integer $n$ such that $A$ has an iterated pullback decomposition with topological dimension $n$. See \cite{Ph} for a detailed exposition of recursive subhomogeneous algebras.
\alten

\altbn{\label{prim}}
Let $A$ be a $C^*$-algebra. We denote by $\Prim \, A$ its primitive ideal space. For $k \in \N$ let $\Prim_k A$ be the subspace of $\Prim \, A$ consisting of kernels of $k$-dimensional irreducible representations. Recall from \cite{Dx}, Propositions 3.6.3 and 3.6.4, that $\bigcup_{k \le n} \Prim_k A$ is closed in $\Prim \, A$, that $\Prim_n A$ is open in $\bigcup_{k \le n} \Prim_k A$ and that $\Prim_n A$ is locally compact Hausdorff.
\alten 

\altbn{\label{characterization}}
There is a nice abstract characterization of recursive subhomogeneous algebras:

\begin{ntheorem}
(\cite{Ph}, Theorem 2.16) For a separable unital $C^*$-algebra $A$ the following are equivalent:\\
(i) $A$ is recursive subhomogeneous of topological dimension not exceeding $n$.\\
(ii) All irreducible representations of $A$ have dimension at most $N$ for some $N \in \N$ and $\dim \Prim_k A \le n$ for $k=1, \ldots, N$.
\end{ntheorem}
\alten

\altbn{\label{rec.subh.}}
We now state the main result of these notes. For the moment we only prove the part which is easy, given the preceding characterization.  

\begin{ntheorem}
Let $A$ be a separable, subhomogeneous $C^*$-algebra. Then we have  
\[
\dr A = \max_k \{\dim \Prim_k A\}\, .
\]
If, additionally, $A$ is unital and $n \in \N$, then $\dr A = n$ iff $A$ is recursive subhomogeneous of topological dimension $n$.\\
\end{ntheorem}

\begin{nproof}
Observe that the first assertion follows from the second: Namely, if $A$ is separable and subhomogeneous, then so is $\tilde{A}$, its smallest unitization. Now by \cite{KW}, Proposition 3.13, we know that $\dr \tilde{A} = \dr A$. Furthermore, it is not hard to see that $\Prim_k \tilde{A} = \Prim_k A$ for $k >1$ and $\Prim_1 \tilde{A}$ is the one-point compactification of $\Prim_1 A$. So, again by \cite{KW}, Proposition 3.13, $\dim \Prim_k A = \dim \Prim_k \tilde{A}$ for all $k$. Thus we only have to show that $\dr \tilde{A} = \max_k \{\dim \Prim_k \tilde{A}\}$. But this follows from the second assertion of the theorem in connection with Theorem \ref{characterization}, since $\tilde{A}$ is separable, subhomogeneous and unital.\\
Therefore, suppose $A$ is separable, unital, subhomogeneous and $\dr A = n$. By \ref{prim},   $\Prim_k A$ is a locally compact Hausdorff spac for any $k \in \N$. Since $A$ is separable, $\Prim_k A$ is second countable and we can find countably many subsets $U_i$ which cover $\Prim_k A$ and each of which has compact closure. It follows from \cite{Fe}, Theorem 3.2, that the corresponding quotients $A_{\overline{U_i}}$ are homogeneous (hence continuous trace) algebras over $\overline{U_i}$. But then \cite{KW}, 3.3 and 3.12, yield the estimate $\dim \overline{U_i} = \dr A_{\overline{U_i}} \le \dr A = n$ for each $i$. Now the countable sum theorem for covering dimension (\cite{HW}, Theorem III.2) says that $\dim \Prim_k A \le n$.\\
The converse turns out to be surprisingly complicated; we postpone the proof to Section 6.  
\end{nproof}
\alten

\altbn{\label{t-crossed-products}}
As an application, we note the following consequence of work by Lin and Phillips; this was our main motivation for studying the decomposition rank of recursive subhomogeneous algebras.
 
\begin{ncor}
Let $M$ be a compact smooth manifold and $\alpha : M \to M$ a minimal diffeomorphism. Then $\cpr (\Ch(M) \times_\alpha \Z) \le \dr (\Ch(M) \times_\alpha \Z) \le \dim M$.
\end{ncor}

\begin{nproof}
By \cite{LP2}, Theorem 1.1, $\Ch(M) \times_\alpha \Z$ can be written as an inductive limit of recursive subhomogeneous algebras of topological dimension at most $\dim M$. The assertion follows from Theorem \ref{rec.subh.} and the permanence properties of the decomposition rank, cf.\ \cite{KW}, 3.3(ii); we have $\cpr A \le \dr A$ for any $A$ by \cite{KW}, Remark 3.2(i).
\end{nproof}
\alten

\altbn
\begin{nremarks}
(i) It would certainly be desirable to find conditions under which we have $\dr \Ch(M) \times_\alpha \Z = \dim M$. However, as for most dimension theories, it is difficult to find lower bounds for the decomposition rank (at least we know from \cite{KW}, 6.1(i), that $\dr A > 0$ unless $A$ is an $AF$ algebra). One way to tackle this problem would be to investigate the $K$-theoretic implications of finite decomposition rank, since the $K$-theory of crossed product $C^*$-algebras as above is quite well understood (cf.\ \cite{LP1}).\\
(ii) It is natural to ask wether similar results still hold for more general crossed products. But then the corollary, as it stands, does not remain true. For example, consider a minimal action on the Cantor set (cf.\ \cite{Ro}, 3.2.12); the crossed product is an $AT$ algebra but it is not $AF$, hence has decomposition rank one (unlike the Cantor set, which is the prototype of a zero-dimensional space). However, an estimate like $\dr (\Ch(X) \times_\alpha \Z) \le \dim X + 1$ is still conceivable for minimal homeomorphisms $\alpha$ of arbitrary compact metrizable spaces $X$. This is tempting not only because of the preceding example, but also because such $C^*$-algebras in some sense might be thought of as  `skew' tensor products of $\Ch(X)$ and $C^*(\Z) = \Ch(S^1)$.
\end{nremarks}
\alten

We close this section with some more examples.
      
\altbn{\label{dimension-drop}}
The unital dimension drop intervals 
\[
\tilde{I}_m := \{ f: [0,1] \to M_m \, | \, f(0), f(1) \in \C \cdot \be_m\}
\]
are recursive subhomogeneous of topological dimension 1 (cf.\ \cite{Ph}), so we have $\dr \tilde{I}_m = 1$.\\
Furthermore, we have $\dr (\Ch(\T) \otimes \tilde{I}_m) = 2$: If $m=1$, then $\Ch(\T) \otimes \tilde{I}_m \cong \Ch( \T \times [0,1])$, and $\T \times [0,1]$ is $2$-dimensional. Otherwise, $\Ch(\T) \otimes \tilde{I}_m$ has only $m$-dimensional and $1$-dimensional irreducible representations and $\Prim_m (\Ch(\T) \otimes \tilde{I}_m) = \T \times (0,1)$ and $\Prim_1 (\Ch(\T) \otimes \tilde{I}_m) = \T \times \{0,1\}$. Therefore, $\max_k \{\Prim_k (\Ch(\T) \otimes \tilde{I}_m)\} = 2$ and the assertion follows from Theorem \ref{rec.subh.}.
\alten

\altbn{\label{dd-extensions}}
Next, consider an extension of the form 
\[
0 \to \Ch_0(\R) \to B \to \Ch(\T) \otimes \tilde{I}_m \to 0
\]
with $B$ unital. By \cite{Ph}, Lemma 2.11, there is a second countable compactification $Y$ of $\R$ with $\dim Y = 1$ and a unital $*$-homomorphism $\pi : \Ch(\T) \otimes \tilde{I}_m \to \Ch(Y \setminus \R)$ such that $B \cong (\Ch(\T) \otimes \tilde{I}_m) \oplus_{\pi, Y \setminus \R} \Ch(Y)$. But then by \cite{Ph}, Proposition 3.2, $B$ is recursive subhomogeneous of topological dimension $2 = \dr B$.
\alten

\altbn{\label{range}}
Let $(G_0,(G_0)_+)$ be a countable ordered abelian group which is weakly unperforated (i.e., if $nx >0$ for some $n \in \N$, then $x>0$), and let $G_1$ be another countable abelian group.\\
It follows from work of Elliott and Thomsen, that there is a (unital) inductive system $B_0 \to B_1 \to \ldots \, $ of $C^*$-algebras as in \ref{dd-extensions}, s.t.\ $A := \lim_\to B_i$ is a simple, separable, unital, nuclear $C^*$-algebra with $(K_0(A), K_0(A)_+,K_1(A)) = (G_0, (G_0)_+, G_1)$; cf.\ \cite{Ro}, Theorem 3.4.4 and Example 3.4.8.\\
In particular, if the Elliott conjecture (see \cite{Ro}, Conjecture 2.2.5) was true, this would imply that $\dr A \le 2$ for all stably finite, simple, separable, unital, nuclear $C^*$-algebras with weakly unperforated $K_0$-groups. 
\alten

\altbn
Recursive subhomogeneous algebras directly generalize the noncommutative $CW$-complexes of \cite{Ped}. A noncommutative $CW$-complex $A$ arises as in Definition \ref{d-rec.subh.}, with the restriction that (at each step) $\Omega$ is the closed unit ball in $\R^n$ for some $n \in \N$ and $X$ is the $(n-1)$-sphere $S^{n-1}$ in $\R^n$. The topological dimension of $A$ (which equals $\dr A$ by Theorem \ref{rec.subh.}) is the highest number $n$ which occurs in this iterated pullback construction.
\alten

\section{The commutative case: an outline}

Since our proof of Theorem \ref{rec.subh.} is regrettably complicated, it might be helpful to study the commutative case first; this will not only serve as a model for the general case, but it also shows which technical difficulties arise and how to circumvent them.

So let $A$ be a separable commutative recursive subhomogeneous algebra; then $A$ can be written as $\Ch(\Omega_0) \oplus_{\pi,X} \Ch(\Omega)$, where $\Omega_0$ and $X \subset \Omega$ are compact metrizable spaces and $\pi : \Ch(\Omega_0) \to \Ch(X)$ is a unital $*$-homomorphism. We have to show that $\dr A \le \max \{ \dim \Omega, \, \dim \Omega_0 \}$.\\
Note that $\hat{A}$ coincides with the pushout $\Omega_0 \coprod_{\pi, X} \Omega$, so the statement is equivalent to saying $\dim (\Omega_0 \coprod_{\pi, X} \Omega) \le \max \{\dim \Omega, \, \dim \Omega_0\}$ (recall from \cite{KW}, Proposition 3.4, that $\dr A = \dim \hat{A}$ if $A$ is commutative). There are several ways to prove this, for example by using the characterization of covering dimension via maps into spheres (cf.\ \cite{HW}, Theorem VI.4). Below we sketch a direct (although, admittedly less elegant) proof which can be generalized to noncommutative recursive subhomogeneous algebras. 

{\it Step 1.} Let $a_1, \ldots, a_k$ be given. Choose a (sufficiently good) c.p.\ approximation $(\C^s, \psi', \varphi')$ (of $\Ch(\Omega_0)$) for $\beta(a_1), \ldots, \beta(a_k)$ such that $\varphi'$ is $n$-decomposable, where $n:= \max (\dim \Omega, \, \dim \Omega_0)$ and $\beta : A \to \Ch(\Omega_0)$ is the projection map. For convenience we may assume that $\varphi'$ is unital. Using the Tietze extension theorem and functional calculus we can find a small neighborhood $Y' \subset \Omega$ of $X$ and a unital $n$-decomposable map $\hat{\varphi} : \C^s \to \Ch(\Omega_0) \oplus_{\pi, X} \Ch(Y')$ such that $\beta \verk \hat{\varphi}$ is close to $\varphi'$. 

{\it Step 2.} $\Omega$ is normal, so there are open sets $V$, $W \subset \Omega$ and a closed set $Y \subset \Omega$ s.t.\ $X \subset W \subset Y \subset V \subset Y'$.\\
Next construct a finite collection $(U_\lambda)_\Lambda$ of open subsets of $\Omega \setminus X$ as follows:\\
(i) $(U_\lambda)_\Lambda$ is $n$-decomposable as a collection of subsets (cf.\ \cite{KW}, Definition 1.4)\\
(ii) $\Omega \setminus W \subset \bigcup_\Lambda U_\lambda$, $U_\lambda \cap (\Omega \setminus W) \neq \emptyset \; \forall \, \lambda$ and $U_\lambda \subset Y$ whenever $U_\lambda \cap W \neq \emptyset$\\
(iii) for each $\bar{\lambda} \in \Lambda$ there is $t_{\bar{\lambda}} \in U_{\bar{\lambda}}$ s.t.\ $t_{\bar{\lambda}} \notin \bigcup_{\Lambda \setminus \{\bar{\lambda}\}} U_\lambda$\\
(iv) the restricted functions $a_l|_{U_\lambda}$ are almost constant for all $l$ and $\lambda$\\
(v) $U_\lambda \subset Y' \; \forall \, \lambda \in \Lambda' := \{ \lambda \in \Lambda \, | \, U_\lambda \cap Y \neq \emptyset \}$\\
(vi) for each $\lambda \in \Lambda'$ there is $i(\lambda) \in \{1, \ldots, s\}$ such that $\frac{1}{4(n+1)} < \hat{\varphi}(e_{i(\lambda)}) \; \forall  \, t \in U_\lambda$.\\
Set $\Lambda^{(1)}:= \{ \lambda \in \Lambda \, | \, U_\lambda \cap W \neq \emptyset \} $ and $\Lambda^{(2)}:= \Lambda \setminus \Lambda^{(1)}$. Choose positive functions $g_\lambda \in \Ch_0(U_\lambda)$ for each $\lambda \in \Lambda$ s.t.\ $ g:= \sum_\Lambda g_\lambda$ satisfies $0 \le g \le 1$ and $g|_{\Omega\setminus W} \equiv 1$. 

{\it Step 3.} Define $\bar{F} := \bar{F}^{(1)} \oplus \bar{F}^{(2)}$ with $\bar{F}^{(1)} := \C^s$, $\bar{F}^{(2)} := \C^{\Lambda^{(2)}}$ and a c.p.c.\ map $\bar{\psi} : A \to \bar{F}$ by $\psi' \verk \beta \oplus (\bigoplus_{\lambda \in \Lambda^{(2)}} \ev_{t_\lambda})$. The u.c.p.\ map $\bar{\varphi}: \bar{F} \to A$ is given as follows: On $\bar{F}^{(1)}$ define
\[
\textstyle
\bar{\varphi}^{(1)}_i := (\be - g) \cdot \hat{\varphi}_i + \sum_{ \{\lambda \in \Lambda^{(1)} \, | \, i(\lambda) = i\} } g_\lambda \cdot \hat{\varphi}_i \, ,
\]
and on $\bar{F}^{(2)}$ set $\bar{\varphi}^{(2)}(e_\lambda) := g_\lambda$. Now $(\bar{F}, \bar{\psi}, \bar{\varphi})$ is a c.p.\ approximation for $a_1, \ldots, a_k$ which has almost the right properties. In particular, $\bar{\varphi}$ has order (in the topological sense) not exceeding $n$ and the restriction $\bar{\varphi}|_{\C^s}$ even is $n$-decomposable. However, we want $\bar{\varphi}$ to be $n$-decomposable, so we need a modification.

{\it Step 4.} Let $\Delta$ denote the full simplex with vertex set $\{1, \ldots, s \} \stackrel{.}{\cup} \Lambda^{(2)}$, then $\bar{\varphi}$ induces a natural surjection $\tau : \Ch(\Delta) \to C := C^*(\bar{\varphi}(\bar{F}))$ by sending the $j$-th coordinate function to the image under $\bar{\varphi}$ of the $j$-th generator of $\bar{F}$. Let $K$ be the minimal subcomplex of $\Delta$ such that $\tau$ factorizes as $\Ch(\Delta) \to \Ch(K) \to C$ and let $J$ be the subcomplex of $K$ generated by the vertex set $\{1, \ldots, s\}$. By construction, $K$ is $n$-dimensional and the $1$-skeleton of $J$ (regarded as a graph) is $(n+1)$-colourable (see \cite{KW}, Section 1, for the relation between graph colourings and decomposability of c.p.\ maps). But now we can apply what might be called relative barycentric subdivision. This is, we form a subdivision $\Sd_J K$ of $K$ such that $J$ is left fixed and the $1$-skeleton of $\Sd_J K$, again regarded as a graph, is $(n+1)$-colourable. If $\Gamma$ denotes the vertex set of $\Sd_J K$, we may replace $\bar{F}$ by $F := \C^\Gamma$ and $\bar{\varphi}$ by a c.p.c.\ map $\varphi$ which is determined by the coordinate functions of $\Sd_J K$. How to obtain $\psi$ from $\bar{\psi}$ is then quite obvious. The triple $(F, \psi, \varphi)$ now is the desired $n$-decomposable c.p.\ approximation of $A$ which shows that in fact $\dr A \le n$. 

It is this last step that causes a lot of technical difficulties in the general setting, since commutativity of $C$ is essential to define a barycentric subdivision and its coordinate functions. Luckily, we may always assume the image of $\bar{\varphi}$ to be commutative `enough' (as will be made precise in the next two sections), so that we can apply the idea of relative barycentric subdivision (Section 5) to certain maps into recursive subhomogeneous algebras. In the last section we will follow the lines of the above argument to prove the remaining part of Theorem \ref{rec.subh.}.

\section{Central lifting and continuous bundles}

In this section we recall the notion of continuous $M_r$-valued $C^*$-algebra bundles and provide a lifting result for the centers of such bundles. This will be the main ingredient for our approximate lifting theorem of so-called piecewise commuting maps, see Section 4.

\altbn
Given $r \in \N$, let $\Mf_r$ denote the set of all $C^*$-subalgebras of $M_r$. By an $\Mf_r$-bundle over a compact space $\Omega$ we mean a map $\Bh : \Omega \to \Mf_r$; we write $\Bh(\Omega)$ for the $C^*$-algebra of continuous selections or, more precisely, 
\[
\Bh(\Omega) := \{ x \in \Ch(\Omega,M_r) \, | \, x(t) \in \Bh(t) \; \forall \, t \in \Omega \}.
\]
$\Bh(\Omega)$ is a $\Ch(\Omega)$-module and we have $\Ch(\Omega) \cdot \Bh(\Omega) = \Bh(\Omega)$.\\
It will be convenient to define a bundle $\Bh$ to be unital if $\Bh(\Omega)$ is a unital $C^*$-algebra; note that with this definition the unit does not necessarily coincide with that of $\Ch(\Omega) \otimes M_r$. $\Bh$ defines a continuous $C^*$-algebra bundle in the sense of \cite{Dx}, 10.3, if $\Bh(t) = \{x(t) \, | \, x \in \Bh(\Omega) \} \; \forall \, t \in \Omega$; in this case we say $\Bh$ is a continuous $\Mf_r$-bundle over $\Omega$. If $\Ah$ is another continuous $\Mf_r$-bundle over $\Omega$ such that $\Ah(t) \subset \Bh(t) \; \forall \, t$, we say $\Ah$ is a subbundle of $\Bh$, $\Ah \subset \Bh$. 
\alten

\altbn{\label{restriction}}
Recall from \cite{Mi} that $\Bh$ is lower semicontinuous, if for every open subset $\Uh \subset M_r$ the set $\{ t \in \Omega \, | \, \Bh(t) \cap \Uh \neq \emptyset \}$ is open in $\Omega$. It is easy to see that a continuous $\Mf_r$-bundle is lower semicontinuous. Conversely, it follows from the Michael selection principle (cf.\ \cite{Mi}, Theorem 3.2'') that if $\Bh$ is lower semicontinuous, then it is a continuous $\Mf_r$-bundle.\\
Now let $X \subset \Omega$ be a closed subspace. The restriction $\Bh|_X$ obviously is a continuous $\Mf_r$-bundle if $\Bh$ is; in this case Michael selection yields that 
\[
\Bh|_X (X) = \{ x \in \Ch(X, M_r) \, | \, \exists \, \bar{x} \in \Bh(\Omega): \bar{x}(t) = x(t) \; \forall \, t \in X \}
\]
and we write $\Bh(X)$ for $\Bh|_X (X)$. In particular we see that $\Bh(X)$ is a quotient of $\Bh(\Omega)$. \\
A similar reasoning shows that, if $\Ah \subset \Bh|_X$ is a continuous subbundle, then $\Ah$ extends to a continuous $\Mf_r$-bundle on all of $\Omega$ by setting $\Ah(t) := \Bh(t)$ for $t \in \Omega \setminus X$.
\alten

\altbn{\label{fibers}}
Let $B \subset \Ch(\Omega)$ be a $C^*$-subalgebra. Then restriction to fibers yields a map $\Bh : \Omega \to \Mf_r$, given by $\Bh(t):= \{b(t) \in M_r \, | \, b \in B\}$. $\Bh$ clearly is l.s.c.; it is not hard to show that $\Bh(\Omega) = B$ iff $\Ch(\Omega) \cdot B = B$.\\
Consider $\Zh(\Bh(\Omega))$, the center of the $C^*$-algebra $\Bh(\Omega)$. Restriction to fibers yields a continuous $\Mf_r$-bundle, denoted by $\Zh(\Bh)$. Because $\Ch(\Omega) \cdot \Zh(\Bh(\Omega)) = \Zh(\Bh(\Omega))$, we have $\Zh(\Bh(\Omega)) = \Zh(\Bh)(\Omega)$. However, if $X \subset \Omega$ is closed, then $\Zh(\Bh(X))$ ($=\Zh(\Bh|_X)(X)$) need not coincide with $\Zh(\Bh)(X)$, as easy examples show (we certainly have $\Zh(\Bh)(X) \subset \Zh(\Bh|_X)(X)$).\\
The main result of the present section says that this problem can be circumvented by making $\Bh(t)$ smaller for $t \notin X$:
\alten

\altbn{\label{central-extension}}
\begin{nlemma}
Let $\Omega$ be compact and metrizable, $X \subset \Omega$ a closed subspace and $\Bh$ a unital continuous $\Mf_r$-bundle over $\Omega$. Then there is a unital continuous subbundle $\Dh \subset \Bh$ such that $\Dh|_X = \Bh|_X$ and $\Zh(\Dh)(X) = \Zh(\Dh(X))$.
\end{nlemma}

\noindent
Although the result looks plausible, it is not so easy to prove; this is because making a bundle smaller will usually affect its continuity. Unfortunately, we only have a rather complicated argument which will make use of a series of ad hoc constructions.
\alten

\altbn{\label{matrix-units}}
First we need some more notation and two simple observations on sets of matrix units of finite-dimensional $C^*$-algebras. If $F= M_{r_1} \oplus \ldots \oplus M_{r_s}$, we say $\{e_1, \ldots,e_k\} \subset F$ is a set of matrix units for $F$ if $k = \sum_{i=1}^s r_i^2$ and 
\[
\{e_1, \ldots, e_k\} = \{ f_i^{(l,m)} \, | \, i=1, \ldots,s; \, 1 \le l,m \le r_i\}\, ,
\]
where $\{ f_i^{(l,m)} \, | \, 1 \le l,m \le r_i\}$ is a set of matrix units for $M_{r_i}$.

\begin{nprop}
For each $r \in \N$ there is $\alpha > 0$ such that the following holds: If, for $k \le r^2$, $F \subset M_r$ is a $k$-dimensional $C^*$-algebra with a set of matrix units $\{e_1, \ldots, e_k\}$ and $q \in F$ is a projection, then $q \in \Zh(F)$ or $\|[e_m,q]\| \ge \alpha$ for some $m \in \{1, \ldots, k\}$.
\end{nprop}

\begin{nproof}
Suppose first we have some $F \subset M_r$ given. Let $\Uh$ be the set of all sets of matrix units for $F$ and $\Ph$ the set of projections in $F$. With the obvious topologies coming from the norm on $F$ these are compact metrizable spaces. Furthermore one checks that $\Ph_0 := \Ph \cap \Zh(F)$ and $\Ph_1 := \Ph \setminus \Ph_0$ are compact subsets of $\Ph$. Define a function $c : \Ph_1 \times \Uh \to \R_+$ by setting
\[
c(p, \{e_1, \ldots, e_k\}) := \max_{m=1, \ldots,k} \| [e_m,p] \| \, ;
\]
this is obviously continuous, hence takes its minimum $\alpha_F$ on some $(p, \{e_1, \ldots, e_k\}) \in \Ph_1 \times \Uh$. But, since $p \notin \Zh(F)$, we have $\alpha_F > 0$. Note that $\alpha_F$ only depends on the isomorphism class of $F$.\\
Set $\alpha := \min \{\alpha_F \, | \, F \subset M_r \mbox{ is a $C^*$-subalgebra} \}$. We still have $\alpha >0$, since $M_r$ has (up to isomorphism) only finitely many subalgebras. Now, for this $\alpha$, the assertion holds by construction.
\end{nproof}
\alten

\altbn{\label{matrix-lifting}}
Let $\Bh$ be a unital continuous $\Mf_r$-bundle over the compact space $\Omega$ and $t \in \Omega$. Set $F:= \Bh(t)$, then $\id_F$ has strict order zero, so by \cite{KW}, 2.4, and since $\Bh(t)$ is a quotient of $\Bh(\Omega)$, $\id_F$ lifts to a c.p.c.\ order zero map $\varphi: F \to \Bh(\Omega)$. Define $g: \R \to \R$ by $g(x) := \min \{2x,1\}$ and set $h:= g(\varphi(\be_F)) \in \Bh(\Omega)$. Then there is a c.p.c.\ order zero map $\hat{\varphi} : F \to \Bh(\Omega)$, given by $\hat{\varphi}(\,.\,):=h \sigma(\,.\,)$, where $\sigma : F \to \Bh(\Omega)''$ is as in \ref{order-zero}. \\
Set $V := \{x \in \Omega \, | \, \|\varphi(\be_F)(x)\| > \halb\}$ and let $\pi_{\overline{V}} : \Bh(\Omega) \to \Bh(\overline{V})$ be the quotient map. Then one checks that $\pi_{\overline{V}} \circ \hat{\varphi}: F \to \Bh(\overline{V})$ is an order zero map with $\pi_{\overline{V}} \circ \hat{\varphi} (\be_F) = \be_{\Bh(\overline{V})}$, hence $\pi_{\overline{V}} \circ \hat{\varphi}$ is a $*$-homomorphism by \ref{order-zero}. In particular, if $e_1, \ldots,e_k$ is a set of matrix units for $F$, then $\pi_{\overline{V}} \circ \hat{\varphi}(e_j)$, $j=1, \ldots,k$, are matrix units in $\Bh(\overline{V})$ s.t.\ $C^*(\pi_{\overline{V}} \circ \hat{\varphi}(e_j)(x) \, | \, j=1, \ldots,k) \cong F \; \forall \, x \in \overline{V}$.
\alten

\altbn{\label{l-extension-2}}
The following construction will be used repeatedly in the proof of Lemma \ref{l-extension-3}. Let $\Bh$ be a unital continuous $\Mf_r$-bundle over a compact metric space $\Omega$, $K \subset \Omega$ a nonempty closed subspace and $\delta > 0$. For $l \in \{1, \ldots, r^2\}$ define 
\[
K^{(l,\Bh)} := \{t \in K \, | \, \Bh(t) \, \mbox{ has vector space dimension } \, l \}
\]
and
\[
\textstyle K^{(<m, \Bh)} := \bigcup_{l < m} K^{(l,\Bh)} \, .
\]
Similar as for the primitive ideal space (cf.\ \ref{prim}), $K^{(<l,\Bh)}$ is closed in $K$ and $K^{(l,\Bh)}$ is open in $K^{(<l+1,\Bh)}$.\\
Now fix some $k \in \{1,\ldots,r^2\}$. For each $t \in K^{(k,\Bh)}$, $\Bh(t)$ is generated by a set of matrix units $\{e_1(t), \ldots, e_k(t)\}$. These lift to matrix units on a neighborhood of $t$ by \ref{matrix-lifting}. But $K^{(k,\Bh)}$ is $\sigma$-compact; it is then straightforward to construct $V_i \subset_{\rm{open}} K$ for $i \in \N$ as follows:
\begin{itemize}
\item[-] $K^{(k,\Bh)} \subset \bigcup_i V_i$
\item[-] $\overline{V_i} \cap K^{(<k,\Bh)} = \emptyset \; \forall \, i$
\item[-] $(V_j \cap K^{(k,\Bh)}) \setminus V_i \neq \emptyset \; \forall \, i < j$
\item[-] $\diam \overline{V_i} \stackrel{i \to \infty}{\longrightarrow} 0$
\item[-] for each $i$ there is $\bi$ s.t. $\overline{V_j} \cap \overline{V_i} = \emptyset \; \forall \, j > \bi$
\item[-] for each $i$ there are matrix units $e_{i,1}, \ldots, e_{i,k} \in \Bh(\overline{V_i})$ s.t.
\begin{eqnarray*}
& & C^*(e_{i,m} \, | \, m=1, \ldots, k) \cong C^*(e_{i,m}(t) \, | \, m=1, \ldots, k) \; \forall \, t \in \overline{V_i} \, , \\
& & \Bh(t) = C^*(e_{i,m}(t) \, | \, m=1, \ldots, k) \; \forall \, t \in \overline{V_i} \cap K^{(k,\Bh)} \, , \\
& & \|e_{i,m} (t) - e_{i,m} (t') \| < \frac{\delta}{2} \; \forall \, t,  t' \in \overline{V_i} \, .
\end{eqnarray*}
\end{itemize}
Just as in \ref{matrix-lifting} (with $\overline{V_0}$ in place of $t$ and $C^*(e_{i,m}, \, m=1, \ldots, k)$ as $F$), we may extend $e_{0,1}, \ldots, e_{0,k}$ to matrix units (which are again denoted by $e_{0,1}, \ldots, e_{0,k}$) in $\Bh(\overline{K_0})$, where $K_0 \subset_{\rm open} \Omega$ is a neighborhood of $V_0$ s.t.\ $\overline{K_0}$ is compact, $\dist (K_0, \, K^{(<k,\Bh)}) = \dist (K_0 \cap K, \, K^{(<k,\Bh)})$ and $K_0 \cap K = V_0$. We may assume $\|e_{0,m} (t) - e_{0,m} (t') \| < \delta$ for $t,  t' \in \overline{K_0}$.\\
Set
\[
\Dh_0(t) := \left\{ 
\begin{array}{ll}
C^*(e_{0,m}(t) \, | \, m=1, \ldots, k) & \mbox{for } \, t \in \overline{\overline{K_0} \setminus K} \\
\Bh(t) & \mbox{else} \, ,
\end{array}
\right.
\]
then $\Dh_0 \subset \Bh$ is a continuous $\Mf_r$-bundle by \ref{restriction}. As before, lift $e_{1,1}, \ldots, e_{1,k}$ to matrix units in $\Dh_0(\overline{K_1})$, where $K_1 \subset_{\rm open} \Omega$ is a neighborhood of $V_1$ with compact closure and such that $\dist (K_1, \, K^{(<k, \Bh)}) = \dist (V_1, \, K^{(<k,\Bh)})$, $K_1 \cap K = V_1$ and $\|e_{1,m} (t) - e_{1,m} (t') \| < \delta$ for $t,  t' \in \overline{K_1}$.\\
Setting
\[
\Dh_n(t) := \left\{ 
\begin{array}{ll}
C^*(e_{n,m}(t) \, | \, m=1, \ldots, k) & \mbox{for } \, t \in \overline{\overline{K_n} \setminus (K \cup \overline{K_0} \cup \ldots \cup \overline{K_{n-1}})} \\
\Dh_{n-1}(t) & \mbox{else} \, ,
\end{array}
\right.
\]
we may proceed inductively to obtain continuous bundles $(\Dh_i)_\N$ and open subsets $(K_i)_\N$ such that $\dist (K_i, \, K^{(<k,\Bh)}) = \dist (V_i, \, K^{(<k,\Bh)})$, $K_i \cap K = V_i$ and $\|e_{i,m} (t) - e_{i,m} (t') \| < \delta$ for $t, \, t' \in \overline{K_i}$. We may clearly assume $\diam \overline{K_i} \to 0$ and (with a little extra effort) that for each $i$ there is $\bi$ s.t.\ $\overline{K_j} \cap \overline{K_i} = \emptyset \; \forall \, j > \bi$. This relation ensures us that, for each $t \in \Omega$, $(\Dh_i(t))_{i \in \N}$ becomes constant as $i$ goes to infinity; therefore it makes sense to define a bundle $\Dh$ over $\Omega$ by $\Dh(t) := \lim_i \Dh_i(t)$. Set $\hat{K} := K \cup (\bigcup_i \overline{K_i})$. This is closed in $\Omega$, since $K_i \cap K \neq \emptyset \; \forall \, i$ and $\diam \overline{K_i} \to 0$. It is straightforward to check that $\Dh$ is lower semicontinuous and that $\Bh(t) = \Dh(t) \; \forall \, t \in K^{(<k+1, \Bh)}$. We have now constructed a unital continuous subbundle $\Dh$ of $\Bh$ over $\Omega$ and a compact set $\hat{K} \subset \Omega$ with the following properties:\\
(i) $\Dh(t) = \Bh(t)$ for $t \in (\Omega \setminus (\overline{\hat{K} \setminus K}))   \cup K^{(<k+1, \Bh)}$\\
(ii) there is a sequence $(K_i)_\N$ of subsets $K_i \subset \hat{K}$ open in $\Omega$ s.t.
\begin{itemize}
\item[-] $\overline{\hat{K} \setminus K} \cup K^{(k,\Bh)} \subset \bigcup_i \overline{K_i}$ 
\item[-] $K^{(k,\Bh)} \subset \bigcup_i K_i$
\item[-] $\overline{K_i} \cap K^{(< k,\Bh)} = \emptyset \; \forall \, i$
\item[-] $K_i \cap K^{(k,\Bh)} \neq \emptyset \; \forall \, i$
\item[-] $\dist (K_i, K^{(< k,\Bh)}) = \dist (K_i \cap K, K^{(< k,\Bh)})$
\item[-] $\diam \overline{K_i} \stackrel{i \to \infty}{\longrightarrow} 0$
\end{itemize}
(iii) for each $i$ there are matrix units $e_{i,1}, \ldots, e_{i,k} \in \Dh(\overline{K_i})$ s.t.
\begin{itemize}
\item[-] $C^*(e_{i,m} \, | \, m=1, \ldots,k) \cong C^*(e_{i,m}(t) \, | \, m=1, \ldots, k) \; \forall \, t \in \overline{K_i}$
\item[-] $\Dh(t) =  C^*(e_{i,m}(t) \, | \, m=1, \ldots, k) \; \forall \, t \in \overline{\overline{K_i} \setminus (K \cup \overline{K_0} \cup \ldots \cup \overline{K_{i-1}})}$; note that 
\item[ ] $\bigcup_i \overline{\overline{K_i} \setminus (K \cup \overline{K_0} \cup \ldots \cup \overline{K_{i-1}})} = \bigcup_i \overline{\overline{K_i} \setminus K}$
\item[-] $\|e_{i,m} (t) - e_{i,m} (t') \| < \delta$ for $i \in \N, \, m =1, \ldots,k, \, t,  t' \in \overline{K_i}$.
\end{itemize}
\alten

\altbn{\label{l-extension-3}}
\begin{nlemma}
Let $X \subset \Omega$ be compact metric spaces, $\Bh$ a unital continuous $\Mf_r$-bundle over $\Omega$ and $\delta > 0$.\\
Then there is a closed neighborhood $K$ of $X$ and a unital continuous $\Mf_r$-bundle $\Dh \subset \Bh|_K$ with the following properties:\\
\\
(i) $\Dh|_X = \Bh|_X$.\\
\\
(ii) For each $k \in \{1, \ldots, r^2 \}$ there is a sequence $(K_i^{(k)})_{i \in \N}$ of open subsets $K_i^{(k)} \subset K$ s.t.
\begin{itemize}
\item[-] $K^{(k,\Dh)} \subset \bigcup_i K_i^{(k)} \subset \bigcup_i \overline{K_i^{(k)}} \subset K^{(r^2,\Dh)} \cup \ldots \cup K^{(k,\Dh)}$
\item[-] $K_i^{(k)} \cap X^{(k,\Dh)} \neq \emptyset \; \forall \, i$
\item[-] $\dist (K_i^{(k)}, \, X^{(< k,\Dh)}) = \dist (K_i^{(k)} \cap K^{(k,\Dh)}, \, X^{(< k,\Dh)})$
\item[-] $\diam (\overline{K_i^{(k)}}) \stackrel{i \to \infty}{\longrightarrow} 0$.
\end{itemize}
(iii) For each $k$ and $i$ there are matrix units $e_{i,1}^{(k)}, \ldots, e_{i,k}^{(k)} \in \Dh(\overline{K_i^{(k)}})$ s.t.\ for $t, \, t' \in \overline{K_i^{(k)}}$
\begin{itemize}
\item[-] $F_i^{(k)} := C^*(e_{i,m}^{(k)} \, | \, m=1, \ldots,k) \cong C^*(e_{i,m}^{(k)} (t) \, | \, m=1, \ldots, k)$
\item[-] $\|e_{i,m}^{(k)} (t) - e_{i,m}^{(k)} (t') \| < \delta$.
\end{itemize}
\end{nlemma}

\begin{nproof}
Apply \ref{l-extension-2} with $X$ in place of $K$ and $r^2$ as $k$ to obtain a continuous bundle $\Dh^{(r^2)} \subset \Bh$, $\hat{K}^{(r^2)} \supset X$, $(K_i^{(r^2)})_{i \in \N}$ and for each $i \in \N$ matrix units $e_{i,1}^{(r^2)}, \ldots , e_{i,r^2}^{(r^2)}$ with the properties of \ref{l-extension-2}. In the next step take $\hat{K}^{(r^2)}$ as $K$, $\Dh^{(r^2)}$ as $\Bh$ and $r^2 - 1$ as $k$ and apply \ref{l-extension-2}. Then proceed inductively to obtain a continuous bundle $\Dh^{(1)} \subset \Bh$ over $\Omega$ after $r^2$ steps. By construction, $\bigcup_{i,k} K_i^{(k)}$ is an open neighborhood of $X$, so it contains a closed neighborhood $K$ of $X$. \\
Define $\Dh := \Dh^{(1)}|_K$ and restrict $K_i^{(k)}$ and $e_{i,m}^{(k)}$ to $K$, then the assertions (ii) and (iii) hold by construction (cf.\ (ii) and (iii) of \ref{l-extension-2}; we only have to check (i).\\
Let $t \in X$ be given, then $t \in X^{(k,\Bh)}$ for some $k \le r^2$. By \ref{l-extension-2}(i) and the construction of $\Dh^{(r^2)}$, we have $\Bh(t) = \Dh^{(r^2)}(t)$. Inductively we see that $\Bh(t) = \Dh^{(r^2)}(t)= \ldots = \Dh^{(k)}(t)$. Furthermore, we have $X \subset \hat{K}^{(r^2)} \subset \ldots \subset \hat{K}^{(1)}$, and therefore
\[
\textstyle
X^{(k,\Bh)} \subset (\hat{K}^{(k)})^{(k,\Bh)} \stackrel{\ref{l-extension-2}(ii)}{\subset} \bigcup_i K_i^{(k)} \subset \hat{K}^{(k)} \subset \hat{K}^{(l)} \; \forall \, l \le k \, .
\]
Now since $\bigcup_i K_i^{(k)}$ is open in $\Omega$, $t$ has a neighborhood $N$ which is open in $\Omega$ and such that $N \subset \hat{K}^{(l)}$ for  $1 \le l \le k$. But then $\overline{\hat{K}^{(l-1)} \setminus \hat{K}^{(l)}} \cap N = \emptyset$, so $t \in \Omega \setminus (\overline{\hat{K}^{(l-1)} \setminus \hat{K}^{(l)}})$, hence $\Dh^{(l-1)}(t) = \Dh^{(l)}(t)$ for $2 \le l \le k$ by \ref{l-extension-2}(i). Therefore, $\Dh^{(1)}(t) = \ldots = \Dh^{(r^2)}(t) = \Bh(t)$.  
\end{nproof}
\alten

\altbn{\label{r-extension}}
\begin{nremark}
It follows from \ref{l-extension-3}(iii) that $\Dh(t) \cong F_i^{(k)}$ for all $t \in \overline{K_i^{(k)}} \cap K^{(k,\Dh)}$. If $\Bh|_X$ is commutative, then so is $\Dh$, since $K_i^{(k)} \cap X^{(k,\Dh)} \neq \emptyset$.
\end{nremark}
\alten

\begin{nproof} 
(of Lemma \ref{central-extension}) Given $r \in \N$, choose $\alpha > 0$ as in Proposition \ref{matrix-units} and fix a metric $d$ on $\Omega$. With $\delta := \frac{\alpha}{12}$ construct a closed neighborhood $K$ of $X$ and a unital continuous bundle $\Dh \subset \Bh|_K$ as in Lemma \ref{l-extension-3}. By \ref{restriction} we may extend $\Dh$ to a continuous bundle on all of $\Omega$ by setting $\Dh(t) := \Bh(t)$ for $t \notin K$. We have to check that indeed $\Zh(\Dh)(X) = \Zh(\Dh(X))$.\\
Note that we have $\Zh(\Bh|_X) \subset \Bh|_X = \Dh|_X$, hence (again by \ref{restriction}) the bundle $\Eh$ over $\Omega$, defined by
\[
\Eh(t) := \left\{ 
\begin{array}{ll}
\Zh(\Bh|_X)(t) & t \in X \\
\Dh(t) & \mbox{else} \, ,
\end{array}
\right.
\]
is a continuous subbundle of $\Dh$. Therefore, we may oncemore apply Lemma \ref{l-extension-3} (together with Remark \ref{r-extension}), this time to $\Eh$ as $\Bh$, to produce a closed neighborhood $G$ (w.l.o.g.\ $G \subset K$) of $X$ and a unital continuous bundle $\Ah \subset \Dh|_G$ with the following properties:\\
\\
(i) $\Ah|_X = \Eh|_X = \Zh(\Dh|_X)$\\
\\
(ii) for each $k \in \{1, \ldots,r \}$ there is a sequence $(G_i^{(k)})_{i \in \N}$ of open subsets $G_i^{(k)} \subset G$ s.t.\  
\begin{itemize}
\item[-] $G^{(k,\Ah)} \subset \bigcup_i G_i^{(k)} \subset \bigcup_i \overline{G_i^{(k)}} \subset G^{(r,\Ah)} \cup \ldots \cup G^{(k,\Ah)}$
\item[-] $G_i^{(k)} \cap X^{(k,\Ah)} \neq \emptyset \; \forall \, i$
\end{itemize}
(iii) for each $k$ and $i$ there are normalized positive elements $q_{i,1}^{(k)}, \ldots, q_{i,k}^{(k)} \in \Ah(G)$ s.t.\  $q_{i,1}^{(k)}|_{\overline{G_i^{(k)}}}, \ldots, q_{i,k}^{(k)}|_{\overline{G_i^{(k)}}}$ are pairwise orthogonal projections (which means that $\C^k \cong C^*(q_{i,m}^{(k)}(t) \, | \, m=1, \ldots, k) = \Ah(t)$ for each $t \in \overline{G_i^{(k)}} \cap G^{(k,\Ah)}$) and such that  $ \| q_{i,m}^{(k)} (t) - q_{i,m}^{(k)} (t') \| < \delta \, \forall \, t,t' \in \overline{G_i^{(k)}}$.\\
By functional calculus we may even assume $q_{i,m}^{(k)}|_{U_i^{(k)}}, \, m=1, \ldots, k$, to be pairwise orthogonal projections on an open neighborhood $U_i^{(k)}$ of $\overline{G_i^{(k)}}$ and that $\| q_{i,m}^{(k)}(t) - q_{i,m}^{(k)}(t') \| < 2 \delta \, \forall \, t, t' \in U_i^{(k)}$.

Now for each $k$ and $i$ there is an open neighborhood $W_i^{(k)} \subset G_i^{(k)}$ of $G_i^{(k)} \cap X$ s.t.\ $q_{i,1}^{(k)} (t), \ldots, q_{i,k}^{(k)}(t) \in \Zh(\Dh(t)) \, \forall \, t \in W_i^{(k)}$:\\
Otherwise there would be a sequence $(t_j)_\N \subset G_i^{(k)}$ with $\dist (t_j, G_i^{(k)} \cap X) \to 0$ and $m \in \{1, \ldots, k\}$ s.t.\ $q_{i,m}^{(k)}(t_j) \notin \Zh(\Dh(t_j)) \, \forall j$. Since $\overline{G_i^{(k)}}$ is compact, we may assume that $t_j \to \bar{t} \in \overline{G_i^{(k)}} \cap X$. But $q_{i,l}^{(k)} \in \Ah(G) \subset \Dh(G)$, $l= 1, \ldots, k$, hence $\Dh(\bar{t})$ has vector space dimension at least $k$, so $\bar{t} \in K^{(r^2,\Dh)} \cup \ldots \cup K^{(k,\Dh)}$. Now since $K^{(< k,\Dh)}$ is closed, we may assume, after passing to a subsequence, that $(t_j)_\N \subset K^{(k',\Dh)}$ for some $k' \ge k$. We have $K^{(k',\Dh)} \subset \bigcup_i K_i^{(k')}$, so each $t_j$ lies in some $K_{i_j}^{(k')}$. Now by Proposition \ref{matrix-units} and, since $q^{(k)}_{i,m}(t_j) \notin \Zh(\Dh(t_j))$, there must be (again after passing to a subsequence) $e_{i_j,m'}^{k'} \in \Dh(\overline{K_{i_j}^{(k')}})$ s.t.\ 
\[
\| [e_{i_j,m'}^{(k')} (t_j), \, q_{i,m}^{(k)}(t_j)] \| \ge \alpha
\]
for all $j$ and some $m' \in \{1, \ldots, k'\}$ (we have $\Dh(t_j) = C^*(e_{i_j, m'}^{(k')} (t_j) \, | \, m' =1, \ldots, k')$ by Remark \ref{r-extension}).\\
Choose $\bar{t}_j \in K_{i_j}^{(k')} \cap X^{(k',\Dh)}$ (this is possible by \ref{l-extension-3}(ii)), then $\bar{t}_j \to \bar{t}$ as well (we have $\diam K_i^{(k')} \stackrel{i \to \infty}{\longrightarrow} 0$) and we may assume that $t_j$, $\bar{t}_j \in U_i^{(k)} \, \forall \, j$. By construction we know that
\[
\| e_{i_j,m'}^{(k')} (t_j) - e_{i_j,m'}^{(k')} (\bar{t}_j) \| < \delta \; \forall \, j
\]
and
\[
\| q_{i,m}^{(k)} (t_j) - q_{i,m}^{(k)} (\bar{t}_j) \| < 2 \delta \; \forall \, j.
\]
But now by the choice of $\delta$ we obtain
\[
\| [ e_{i_j,m'}^{(k')} (\bar{t}_j), \, q_{i,m}^{(k)} (\bar{t}_j)] \| \ge \frac{\alpha}{2} \, .
\]
On the other hand, $e_{i_j,m'}^{(k')} (\bar{t}_j) \in \Bh(\bar{t}_j)$ and $q_{i,m}^{(k)} (\bar{t}_j) \in \Zh(\Bh|_X)(\bar{t}_j) \subset \Zh(\Bh(\bar{t}_j))$, a contradiction, hence $q_{i,1}^{(k)} (t), \ldots, q_{i,k}^{(k)}(t) \in \Zh(\Dh(t))$ for all $t$ in some open neighborhood $W_i^{(k)} \subset G_i^{(k)}$ of $G_i^{(k)} \cap X$; we may assume that $W_i^{(k)}$ is open in $\Omega$.

Now take an arbitrary $t \in X$, then $t \in W_i^{(k)} \cap X^{(k, \Ah)} \subset G_i^{(k)} \cap G^{(k, \Ah)}$ for some $i$ and $k$. Choosing $h \in \Ch_0(W_i^{(k)}) \subset \Ch(\Omega)$ with $h(t)=1$ we see  from (iii) that $\Ah(t) = C^*(q_{i,m}^{(k)}(t) \, | \, m=1, \ldots,k) = C^*((h \cdot q_{i,m}^{(k)})(t) \, | \, m=1, \ldots,k)$.\\
But $h \cdot q_{i,m}^{(k)} \in \Ah(G) \subset \Dh(G) \; \forall \, m$, and since $(h \cdot q_{i,m}^{(k)}) (t') \in \Zh(\Dh(t')) \, \forall \, t' \in \Omega$, we have that $h \cdot q_{i,m}^{(k)} \in \Zh(\Dh)(G)$, $m =1, \ldots, k$. In particular we obtain  
\[
\Zh(\Dh|_X)(t) = \Ah(t) = C^*((h \cdot q_{i,m}^{(k)})(t), \, m=1, \ldots,k) \subset \Zh(\Dh)(t) \; \forall \, t \in X \, ;
\]
it follows that $\Zh(\Dh|_X)(X) \subset \Zh(\Dh)(X)$. The inclusion $\Zh(\Dh)(X) \subset \Zh(\Dh|_X)(X)$ is trivial; hence the proof is complete.
\end{nproof}

\section{Piecewise commuting maps}

\altbn{\label{d-piecewise-commuting}}
\begin{ndefn}
Let $A, \, F$ be $C^*$-algebras, $F = M_{r_1} \oplus \ldots \oplus M_{r_s}$, and $\varphi : F \to A$ a c.p.\ map. If there is an order $\prec$ on $\{1, \ldots, s\}$ s.t.
\[
[ \varphi(\be_{M_{r_j}}), \, \varphi(M_{r_i})] =0
\]
for $i \prec j$, we say $\varphi$ is piecewise commuting (p.c.) with respect to $\prec$. 
\end{ndefn}

\noindent
The reason for introducing this somewhat technical concept is, that it admits use of Lemma \ref{l-piecewise-commuting} (which yields a decomposition of $\varphi$ into pairwise orthogonal parts) and an approximate lifting result (Proposition \ref{pc-extension}) at the same time. 
\alten

\altbn{\label{l-piecewise-commuting}}
The following will play an important r\^ole in the proof of Theorem \ref{rec.subh.}. By $\chi_\mu$ we denote the characteristic function of $[\mu, \infty)$.
 
\begin{nlemma}
Let $A$ be a unital $C^*$-algebra, $F = M_{r_1} \oplus \ldots \oplus M_{r_s}$ and $\varphi: F \to A$ a c.p.\ and p.c.\ map. Suppose there is $\mu > 0$ s.t.\ the projections $q_i := \chi_\mu (\varphi(\be_{M_{r_i}})) \in C^*(\varphi(\be_{M_{r_i}})) \subset A$ exist for all $i$ and $\sum_{i=1}^s q_i \in A$ is invertible. \\
Then there are pairwise orthogonal projections $p(i) \in C^*(q_j \, | \, j \in \{1, \ldots, s\})$, $i= 1, \ldots, s$, s.t.\
\[
{\textstyle \sum_{j=1}^s} \, p(j) = \be_A , \; [p(i), \, \varphi(F)] =0 \mbox{ and } p(i) q_i = p(i)\, .
\]
\end{nlemma}

\begin{nproof}
Set $\Sigma^{(0)} := \emptyset$, $B^{(0)} := C^*( q_j \, | \, j \in \Sigma^{(0)}) (= \{0\})$ and $h^{(0)} := \sum_{i \in \Sigma^{(0)}} p(i) = \be_{B^{(0)}} = 0$. Suppose that for $m = 0, \ldots, l$ we have constructed pairwise disjoint $\Sigma^{(m)} \subset \{1, \ldots,s\}$ and pairwise orthogonal projections $p(i) \in C^*(q_j \, | \, j=1, \ldots, s)$, $i \in \Sigma^{(m)}$, s.t.\ $h^{(m)} := \sum_{i \in \Sigma^{(0)} \cup \ldots \cup \Sigma^{(m)}} p(i) = \be_{B^{(m)}}$, where $B^{(m)} := C^*(q_j \, | \, j \in \Sigma^{(0)} \cup \ldots \cup \Sigma^{(m)})$, and such that $[ p(i), \, \varphi(F)] =0 $ and $p(i) \le q_i$. \\
If $h^{(l)} = \be_A$, we are done. Otherwise, define $\varphi^{(l+1)} : F \to (\be -h^{(l)}) A (\be - h^{(l)})$ by $\varphi^{(l+1)}(\, . \,) := (\be - h^{(l)}) \varphi(\, . \, )$, then $\varphi^{(l+1)}$ is c.p.\ and p.c., because $\varphi$ is and $h^{(l)}$ commutes with $\varphi(F)$. Let 
\[
\Sigma^{(l+1)} := \{ i \in \{1, \ldots,s\} \setminus (\Sigma^{(0)} \cup \ldots \cup \Sigma^{(l)}) \, | \, (\be - h^{(l)}) q_i \in (\varphi^{(l+1)}(F))' \} \, ;
\]
using that $\varphi^{(l+1)}$ is p.c.\ and that $\sum_{i \notin (\Sigma^{(0)} \cup \ldots \cup \Sigma^{(l)})} (\be - h^{(l)}) q_i $ is invertible in $(\be - h^{(l)}) A (\be - h^{(l)})$, one checks that $\Sigma^{(l+1)} \neq \emptyset$.\\
$C^*((\be-h^{(l)})q_j \, | \, j \in \Sigma^{(l+1)})$ is a finite-dimensional abelian $C^*$-algebra, so it contains (not necessarily nonzero) pairwise orthogonal projections $p(i)$, $i \in \Sigma^{(l+1)}$, s.t.\ $p(i) \le (\be - h^{(l)}) q_i$ and $\sum_{j \in \Sigma^{(l+1)}} p(j)$ is the unit of $C^*((\be-h^{(l)})q_j \, | \, j \in \Sigma^{(l+1)})$. Because $\varphi(F) = h^{(l)} \varphi(F) + (\be -h^{(l)}) \varphi(F)$, we have $p(i) \in \varphi(F)' \; \forall \, i \in \Sigma^{(l+1)}$.

Thus we may proceed inductively with the construction of $\Sigma^{(m)}$ and $p(i)$; but $\sum_{i \in \Sigma^{(m)}} p(i) \neq 0$ for $m > 0$ (unless $h^{(m-1)} = \be_A$), so there is $k$ s.t. 
\[
\textstyle
\sum_{l=0}^k \sum_{i \in \Sigma^{(l)}} p(i) = \be_A \, .
\]
Set $p(i) := 0$ for $i \notin \Sigma^{(0)} \cup \ldots \cup \Sigma^{(k)}$, then we are done.
\end{nproof}
\alten

\altbn{\label{p.c.weakly-stable}}
In \cite{KW}, Proposition 2.6, it was pointed out that $n$-decomposable maps are weakly stable. We adjust this result to the case of p.c.\ and $n$-decomposable maps.

\begin{nprop}
For any finite-dimensional $C^*$-algebra $F = M_{r_1} \oplus \ldots \oplus M_{r_s}$ and $\varepsilon > 0$ there is $\gamma > 0$ such that the following holds:\\
Let $B$ be a $C^*$-algebra, $\varphi : F \to B$ c.p.c. and p.c.\ such that, for some decomposition $\coprod_{k=1, \ldots, n} I_k$ of $\{1, \ldots, s\}$, $\|\varphi(\be_i) \varphi(\be_j)\| < \gamma$ whenever $i, \, j \in I_k$ for some $k$. Furthermore, assume that $\varphi$ has strict order zero on each summand of $F$.\\
Then there is a c.p.c.\ map $\hat{\varphi} : F \to B$ which is p.c.\ and $n$-decomposable (w.r.t.\  the given decomposition) and s.t.\ $\|\varphi(x) - \hat{\varphi}(x)\| \le \varepsilon \|x\| \; \forall \, x \in F_+$.
\end{nprop}

\begin{nproof}
Set $\gamma := \frac{\varepsilon^2}{s^2}$ and $\hat{h}_i := (\varphi_i(\be_i) - \gamma^\halb)_+$, $i=1, \ldots,s$. Then by \ref{order-zero} for each $i$ there is a map $\hat{\varphi}_i :M_{r_i} \to B$ s.t.\ $\ord \hat{\varphi}_i = 0$, $\|\hat{\varphi}_i - \varphi_i\| \le \gamma^\halb$ and $\hat{\varphi}_i = \hat{h}_i$. This defines a c.p.c.\ map $\hat{\varphi} : F \to B$, which again is p.c.\ and satisfies $\|\hat{\varphi}(x)-\varphi(x)\| \le s \cdot \gamma^\halb\|x\| = \varepsilon \|x\| \; \forall \, x \in F_+$. Furthermore, if $\hat{\varphi}_i (\be_i) \hat{\varphi}_j (\be_j) \neq 0$ for some $i,j \in \{1, \ldots,s\}$, there is a character $\rho$ of the commutative $C^*$-algebra $C^*(\varphi(\be_i), \varphi(\be_j)) \subset B$ s.t.\ $\rho(\hat{h}_i) \rho(\hat{h}_j) \neq 0$. But then $\rho(\varphi(\be_i)), \rho(\varphi(\be_j)) \ge \gamma^\halb$, hence $\|\varphi(\be_i) \varphi(\be_j)\| \ge \gamma$. It follows that $\hat{\varphi}$ is $n$-decomposable with respect to the given decomposition $\{1, \ldots,s\} = \coprod_{k=1, \ldots,s} I_k$. 
\end{nproof}
\alten

\altbn{\label{commutative-perturbation}}
For later use we note the following simple observation:

\begin{nprop}
Let $A$ be a commutative $C^*$-algebra and $a, h \in A$ positive normalized elements satisfying $\| a - h \| < \delta$ for some $0 < \delta \le 1$. Then there are $h' \in (C^*(h))_+$ and $g \in (C^*(h,  a))_+$ with $\| h - h' \| < \delta $ and $\| g \| \le 1 $ s.t.\ $g a = h'$.
\end{nprop}

\begin{nproof}
We may assume $A= C^*(h,a) = \Ch_0(X)$ for some locally compact space $X$. Set $h':=(h-\delta)_+$, then $\|h'-h\|\le \delta$ and $h' \le a$. Define a continuous function $f$ on $\R_+$ by 
\[
f(t) := \left\{ 
\begin{array}{ll}
0, & t=0 \\
t^{-1}, & t\ge \delta - \|a-h\| \\
\mbox{linear,} & \mbox{else}
\end{array}
\right.
\]
and set $g:= f(a) h'$. Then $0 \le g \le f(a) a \le \be$ and one checks that $ga = h'$.
\end{nproof}
\alten

\altbn{\label{pc-extension}}
It is the following approximate lifting result for which we made all the effort in Section 3.

\begin{nprop}
Let $A$ and $B$ be separable recursive subhomogeneous $C^*$-algebras s.t.\ $A$ is of the form $B \oplus_{\pi, X} (\Ch(\Omega) \otimes M_r)$ for $X \subset \Omega$ compact and metrizable spaces and $\pi : B \to \Ch(X) \otimes M_r$ a unital $*$-homomorphism. Suppose that $F$ is a finite-dimensional $C^*$-algebra, that $\varphi : F \to B$ is c.p.c., p.c.\ and $n$-decomposable and that $\be_B \in C^*(\varphi(F))$.\\
Then for any $\alpha > 0$ there are a closed neighborhood $Y$ of $X$ and a c.p.c.\ map $\hat{\varphi} : F \to B \oplus_{\pi, X} (\Ch(Y) \otimes M_r)$, again p.c.\ and $n$-decomposable, s.t.
\[
\| \varphi(x) - \beta \hat{\varphi}(x) \| < \alpha \|x\| \; \forall \, x \in F_+\, ,
\]
where $\beta$ is the projection map onto $B$. 
\end{nprop}

\begin{nproof}
We have $B \subset \Ch(\Omega_0) \otimes M_{r_0}$ for some compact metrizable $\Omega_0$ and $r_0 \in \N$. Set $r' := \max\{r_0,r\}$, $\Omega' := \Omega \stackrel{.}{\cup} \Omega_0$ and $X' := X \stackrel{.}{\cup} \Omega_0$; then, using the upper left corner embeddings of $M_{r_0}$ and $M_r$ into $M_{r'}$, we obtain an injection $A \hookrightarrow B \oplus (\Ch(\Omega) \otimes M_r) \hookrightarrow \Ch(\Omega') \otimes M_{r'}$. Restriction to fibers as in \ref{fibers} now defines a unital continuous $\Mf_{r'}$-bundle $\Ah$ over $\Omega'$.

We may assume $F = M_{r_1} \oplus \ldots \oplus M_{r_s}$ and $[\varphi(\be_j), \varphi(M_{r_i})] = 0$ if $i \le j$. Define a unital continuous $\Mf_{r'}$-bundle $\Eh_s \subset \Ah$ over $\Omega'$ by setting 
\[
\Eh_s(t) := \left\{ 
\begin{array}{ll}
C^*(\varphi({\textstyle \bigoplus_{i=1}^s M_{r_i}})(t), \be_{\Ah(t)})\, , & t \in X' \\
\Ah(t) (= M_r) \, , & t \in \Omega' \setminus X' \, .
\end{array}
\right.
\]
Now apply Lemma \ref{central-extension} (with $\Eh_s$ as $\Bh$, $\Omega'$ as $\Omega$ and $X'$ as $X$) to obtain a unital continuous subbundle $\Dh_s \subset \Eh_s$ s.t.\ $\Dh_s|_{X'} = \Eh_s|_{X'}$ and $\Zh(\Dh_s)(X') = \Zh(\Dh_s(X')) = \Zh(\Eh_s(X'))$.\\
If $\Eh_j$ and $\Dh_j$ are defined, we obtain $\Eh_{j-1}$ and $\Dh_{j-1}$ as follows:\\
Set
\[
\Eh_{j-1}(t) := \left\{ 
\begin{array}{ll}
C^*(\varphi({\textstyle \bigoplus_{i=1}^{j-1} M_{r_i}})(t), \be_{\Ah(t)})\, , & t \in X' \\
\Dh_j(t)  \, , & t \in \Omega' \setminus X' 
\end{array}
\right.
\]
and apply \ref{central-extension}, this time with $\Eh_{j-1}$ as $\Bh$, to obtain $\Dh_{j-1} \subset \Eh_{j-1} \subset \Dh_j$ s.t.\ $\Dh_{j-1}|_{X'} = \Eh_{j-1}|_{X'}$ and $\Zh(\Dh_{j-1})(X') = \Zh(\Dh_{j-1}(X')) = \Zh(\Eh_{j-1}(X'))$. Proceed inductively to construct unital continuous bundles $\Dh_1 \subset \ldots \subset \Dh_s \subset \Ah$ over $\Omega'$ satisfying $\Zh(\Dh_j)(X') = \Zh(\Eh_j(X'))$ for all $j$.

By \ref{fibers} we have $\Eh_j(X') = \Ch(X') \cdot C^*( \varphi({\textstyle \bigoplus_{i=1}^j M_{r_i}}), \be_{\Ah(X')})$, hence, by our assumption on $\varphi$, $\varphi(\be_j) \in \Zh(\Eh_j(X'))$. But $\Zh(\Eh_j(X')) = \Zh(\Dh_j)(X')$ is a quotient of $\Zh(\Dh_j)(\Omega')$, so $\varphi(\be_j)$ lifts to some $h_j \in \Zh(\Dh_j)(\Omega')$ with $0 \le h_j \le \be$. Note that $h_j$ commutes with $\Dh_i(\Omega')$ whenever $i \le j$.\\
Furthermore, by \cite{Wi2}, Proposition 1.2.4, for each $j$ there is a c.p.c.\ lift $\bar{\varphi}_j : M_{r_j} \to \Dh_j(\Omega')$ of $\varphi_j$ s.t.\ $\ord \bar{\varphi}_j = 0$.

Apply Proposition \ref{p.c.weakly-stable} to obtain (from $F$ and $\varepsilon := \frac{\alpha}{2}$) a $\gamma > 0$ (we may assume $\gamma \le \frac{\alpha}{2}$. Choose a compact neighborhood $Y' \subset \Omega$ of $X'$ s.t.\ 
\[
\| \bar{\varphi}_j(\be_j)|_{Y'} - h_j|_{Y'} \| < \min \{{\textstyle \frac{\gamma}{4s}},1\} \; \forall \, j \in 1, \ldots,s \, ,
\]

Now by Proposition \ref{commutative-perturbation} there are $0 \le g_j \in C^*(h_j|_{Y'}, \bar{\varphi}_j(\be_j)|_{Y'})$ and $0 \le h_j' \in C^*(h_j|_{Y'})$ as follows:
\begin{itemize}
\item[-] $\|g_j\| \le 1$
\item[-] $\| h_j'(t) - h_j(t) \| \le \frac{\gamma}{4s} \; \forall \, t \in Y'$
\item[-] $g_j(t) \bar{\varphi}_j(\be_j)(t) = h_j'(t) \; \forall \, t \in Y'$.
\end{itemize}
Note that $h_j' \in \Zh(\Dh_j)(Y')$ for all $j$,  and that $[ g_j(t), \, \bar{\varphi}_j(M_{r_j})(t)] = 0 \; \forall \, t \, \in Y'$.\\
Define $\tilde{\varphi}_j : M_{r_j} \to \Dh_j(Y')$ by
\[
\tilde{\varphi}_j (x)(t) := g_j(t) \cdot \bar{\varphi}_j (x)(t) \mbox{ for } t \in Y' \, ,
\]
then $\ord \tilde{\varphi}_j = 0$ and $\tilde{\varphi}_j(\be_j) = h_j'$. Since $\ord \bar{\varphi}_j = 0$, by \ref{order-zero} there exist $*$-homomorphisms $\sigma_j : M_{r_j} \to D_j(Y_1)''$ s.t.\ $\bar{\varphi}_j (x)(t) = \bar{\varphi}_j (\be_j)(t) \cdot \sigma_j(x)(t)$, hence $\tilde{\varphi}_j(x)(t) = h_j'(t) \cdot \sigma_j(x)(t)$. As a consequence,
\[
\| \tilde{\varphi}_j (x)(t) - \bar{\varphi}_j(x)(t) \| \le \|h'_j(t) - \bar{\varphi}_j(\be_j)(t)\|  < \frac{2 \gamma}{4s} \; \forall \, t \in Y', \, x \in M_{r_j} \mbox{ with } \| x \| \le 1 \, .
\]
We now have a c.p.\ and p.c.\ map $\tilde{\varphi} : F \to \Dh(Y')$ which has strict order zero on each summand of $F$ and satisfies $\| \tilde{\varphi}(x)(t) - \bar{\varphi}(x)(t) \| < s \cdot \frac{\gamma}{2s} = \frac{\gamma}{2}$ for $x \in F$ with $\|x\| \le 1$ and $t \in Y'$. It is not hard to see that, by making $Y'$ smaller if necessary, we may assume that $\|\tilde{\varphi}(\be_i) \tilde{\varphi}(\be_j)\| < \gamma$ whenever $\varphi(\be_i) \varphi(\be_j) = 0$. Thus by Proposition \ref{p.c.weakly-stable} there is a c.p.c.\ and p.c.\ map $\hat{\varphi} : F \to \Dh(Y')$ which is $n$-decomposable and satisfies
\[
\| \hat{\varphi}(x) - \tilde{\varphi}(x) \| \le \frac{\alpha}{2} \|x\| \; \forall \, x \in F_+ \, .
\]
Finally we see that, for all $x \in F_+$ with $\|x\| \le 1$ and $t \in X$,
\[
\| \hat{\varphi}(x)(t) - \varphi(x)(t) \| \le \frac{\gamma}{2} + \frac{\alpha}{2} < \alpha \, . 
\]
\end{nproof}
\alten

\section{Relative barycentric subdivision}

Another reason for introducing p.c.\ maps is that they are accessible to a technique which might be called relative barycentric subdivision. In the proof of Theorem \ref{rec.subh.} it will be used to obtain an $n$-decomposable p.c.\ map from a p.c.\ map satisfying a certain order condition. First, we need some notation. See \cite{Sp}, Chapter 3, for an introduction to abstract simplicial complexes.

\altbn
Following \cite{Sp}, Section 3.1, by a finite simplicial complex $K$ we mean a collection of subsets of a finite vertex set $V(K)$ such that, if $\kappa \in K$, then $\kappa' \subset \kappa$ implies $\kappa' \in K$ and such that, if $\nu \in V(K)$, then $\{\nu\} \in K$. We say $\kappa \in K$ is  $k$-face, if its cardinality is $k+1$. Note that the map $\nu \mapsto \{\nu\}$ is a bijection between $V(K)$ and the $0$-faces of $K$.\\
Let $|K|$ be the geometric realization of $K$; by definition, $|K|$ consists of those points $t \in [0,1]^{V(K)}$ for which ${\textstyle \sum_{\nu \in V(K)}} t_\nu =1$ and for which the sets $\{\nu \in V(K) \, | \, t_\nu \neq 0 \}$ are faces of $K$. Then $|K|$ is a subset of the standard simplex $\{t \in [0,1]^{V(K)} \, | \, \sum_{\nu \in V(K)} t_\nu = 1\}$ in $\R^{V(K)}$ (which in turn is the geometric realization $|\Delta^{V(K)}|$ of the full simplex with vertex set $V(K)$, $\Delta^{V(K)}$).\\
There is a canonical open covering $(A_\nu)_{\nu \in V(K)}$ of $|K|$ which comes from open stars around vertices in the standard simplex in $\R^{V(K)}$. More precisely, 
\[
A_\nu =\{t \in |K|\, | \, t_\nu \neq 0\} \, .
\]
If $K'$ is a simplicial complex such that $V(K') \subset V(K)$ and $K' \subset K$, we say $K'$ is a subcomplex of $K$ and identify $|K'|$ with the corresponding subspace of $|K|$. In particular, each $\kappa \in K$ defines (the geometric realization of) a face $|\kappa| \subset |K|$.
\alten

\altbn
Let  $L$ be another simplicial complex, and let $\tau : V(K) \to |L|$ be a map such that, for all $\kappa \in K$, the convex combinations $\{ \sum_{\nu \in \kappa} \lambda_\nu \tau(\nu) \, | \, \lambda_\nu \ge0, \, \sum \lambda_\nu = 1\}$ lie in $|L|$.  Then $\tau$ induces a map $\bar{\tau} :|K| \to |L|$ (called linear) by
\[ 
\bar{\tau}(t) = {\textstyle \sum_{\nu \in V(K)}} t_\nu \tau (\nu)\, .
\]
Denote the coordinate functions of $|K|$ by $\bar{h}_\nu$, $\nu \in V(K)$; by the Stone--Weierstra{\ss} Theorem, these generate $\Ch(|K|)$ as a $C^*$-algebra. Furthermore, they induce a natural u.c.p.\ map $\bar{h}: \C^{V(K)} \to \Ch(|K|)$.
\alten

\altbn{\label{subdivision}}
Let $\Sigma = \Sigma^{(1)} \stackrel{.}{\cup} \Sigma^{(2)}$ be finite sets and $\Sigma^+ := \Sigma \stackrel{.}{\cup} \{*\}$ the disjoint union of $\Sigma$ with a single point, then $\Delta^{\Sigma^+}$ contains $\Delta^{\Sigma^{(1)}}$ and $\Delta^{\Sigma^{(2)}}$ in the obvious way (again, by $\Delta^M$ we mean the full simplex with vertex set $M$).\\
Suppose $C$ is a unital commutative $C^*$-algebra generated by positive elements $h_\sigma$, $\sigma \in \Sigma^+$, s.t.\ $\sum_{\Sigma^+} h_\sigma = \be_C$; again we regard $h$ as a u.c.p.\ map: $\C^{\Sigma^+} \to C$.\\
Let $K$ consist of the one-point subsets of $\Sigma^+$ and of those subsets $\{\sigma_0, \ldots, \sigma_l\}$ of $\Sigma^+$ for which $h_{\sigma_0} \ldots h_{\sigma_l} \neq 0$.  This defines a simplicial complex with vertex set $V(K) = \Sigma^+$. Note that we have a canonical surjection $\Ch(|K|) \to C$ and that $h$ factorizes as $\C^{\Sigma^+} \stackrel{\bar{h}}{\to} \Ch(|K|) \to C$.

\begin{nprop}
Let $C$, $h$, $K$, $\Sigma^+ (=V(K))$ and $\bar{h}$ be as above; let $J:= K \cap \Delta^{\Sigma^{(1)}} \subset K$ be the subcomplex of $K$ generated by $\Sigma^{(1)}$. Suppose that the generators of $C$ satisfy $h_{\sigma_0} \ldots h_{\sigma_{n+1}} = 0$ for any choice of distinct elements $\sigma_0, \ldots, \sigma_{n+1} \in \Sigma$ and that $h|_{\C^{\Sigma^{(1)}}}$ is $n$-decomposable.\\
Then there is a simplicial complex $\Sd_J K$ with vertex set $\Gamma := V(\Sd_J K)$ satisfying the following conditions:\\
(i) $\Sigma^+ \subset \Gamma \subset \{\mbox{faces of } K \}$; if $\gamma \in \Gamma \setminus \Sigma^+$, then $\gamma$ (as a face of $K$) intersects $\Sigma^{(2)} \subset V(K)$.\\
(ii) The map $\beta: \Gamma \to |K|$, sending each $\gamma \in \Gamma$ to the barycenter $\sum_{\nu \in \gamma} \frac{1}{\card \gamma} \cdot |\nu| \in |K|$ of the corresponding face in $|K|$, induces a linear homeomorphism $\bar{\beta} : |\Sd_J K| \to |K|$.\\
(iii) $\Sd_J K \cap \Delta^{\Sigma^{(1)}} = J$ (so, regarding $|J|$ as a subspace of $|\Sd_J K|$, we have $\bar{\beta}|_{|J|} = \id_{|J|}$).\\
(iv) The coordinate functions $\bar{k}_\gamma \in \Ch(|\Sd_J K|)$, $\gamma \in \Gamma$, induce a u.c.p.\ map $\bar{k} : \C^\Gamma \to \Ch(|\Sd_J K|)$ s.t.\ $\bar{k}|_{\C^{\Gamma \setminus \{*\}}}$ is $n$-decomposable and s.t.\ $\bar{k}_\gamma (t) = \bar{h}_\gamma \verk \bar{\beta}(t)$ for all $\gamma \in \Sigma^{(1)}$ and $t \in |J| \subset |\Sd_J K|$. Moreover, $\sum_{\Gamma \setminus \{*\}} \bar{k}_\gamma = \sum_{\sigma \in \Sigma} \bar{h}_\sigma \verk \bar{\beta}$.  
\end{nprop}
\alten

\altbn{\label{properties-subdivision}}
Before proving the proposition, we derive some consequences which will be needed later on. First, consider the u.c.p.\ map 
\[
k : \C^\Gamma \stackrel{\bar{k}}{\to} \Ch(|\Sd_J K|) \stackrel{\bar{\beta}_*}{\to} \Ch(|K|) \to C \; (e_\gamma \mapsto k_\gamma \in C) \,;
\]
the restriction of $k$ to $\C^{\Gamma \setminus \{*\}}$ is a composition of an $n$-decomposable map and a $*$-homomorphism, hence again $n$-decomposable.

Because of \ref{subdivision}(i) and (ii),  we can choose a function $\nu : \Gamma \setminus ({\Sigma^{(1)}})^+  \to \Sigma^{(2)}$ with $\nu(\sigma) = \sigma $ if $\sigma \in \Sigma^{(2)}$ and s.t.\ $\nu(\gamma)$ is a vertex of $\gamma \; \forall \, \gamma \in \Gamma \setminus \Sigma^+$. This in particular means that $\bar{\beta}$ maps the open star around $|\gamma|$ in $|\Sd_J K|$ to the open star around $|\nu(\gamma)|$ in $|K|$ for all $\gamma \in \Gamma \setminus (\Sigma^{(1)})^+$.\\
Define 
\[
\Gamma' := \{ \gamma \in \Gamma \setminus (\Sigma^{(1)})^+ \, | \, \bar{k}_\gamma \cdot (\bar{h}_\nu \verk \bar{\beta} \mbox{ for some } \nu \in \Sigma^{(1)} \}
\]
and $\Gamma'' := \Gamma \setminus ((\Sigma^{(1)})^+ \cup \Gamma')$, then $\Gamma = \{*\} \stackrel{.}{\cup} \Sigma^{(1)} \stackrel{.}{\cup} \Gamma' \stackrel{.}{\cup} \Gamma''$.\\
For our function $\nu : \Gamma' \cup \Gamma'' \to \Sigma^{(2)}$ we see that $k_\gamma \in \Jh(h_{\nu(\gamma)}) \subset C \; \forall \, \gamma \in \Gamma' \cup \Gamma''$ ($\Jh(M) \subset C$ denotes the ideal generated by $M \subset C$).\\
Using \ref{subdivision}(ii), (iii), (iv) and the Stone--Weierstra{\ss} Theorem one checks that
\[
k_\gamma \in \left\{ 
\begin{array}{ll}
C^*(h_\gamma, \,  h_{\gamma'}, \, h_{\gamma'}h_{\gamma''}, \, \be_C  \, | \, \gamma'  \in \Sigma^{(2)},\, \gamma'' \in \Sigma^{(1)} \}) \cap \Jh(h_\gamma) & \mbox{if } \gamma \in \Sigma^{(1)} \\
\Jh(h_{\nu(\gamma)}) & \mbox{if } \gamma \in \Gamma' \\
C^*(  h_\sigma, \, \be_C \, | \, \sigma \in \Sigma^{(2)}  ) \cap \Jh(h_{\nu(\gamma)}) & \mbox{if } \gamma \in \Gamma'' \, . 
\end{array} \right.
\]
For the obvious u.c.p.\ map
\[
\varrho : \C^{\Sigma^+} \stackrel{\bar{h}}{\to} \Ch(|K|) \stackrel{\bar{\beta}_{*}}{\to} \Ch(|\Sd_J K|) \stackrel{\ev_\Gamma}{\to} \C^\Gamma
\]
we have $\bar{k} \verk \varrho = \bar{h}$, hence $k \verk \varrho = h$, as follows from the linearity of $\bar{\beta}$. Furthermore, we have $\varrho(\C^{\Sigma}) \subset \C^{\Gamma \setminus \{*\}}$. Finally, by (iii) $\id_{\C^{\Sigma^{(1)}}}$ factorizes as
\[
\C^{\Sigma^{(1)}} \hookrightarrow \C^{\Sigma^+} \stackrel{\varrho}{\to} \C^\Gamma \to \C^{\Sigma^{(1)}}
\]
and
\[
\varrho^{(3,2)} : \C^{\Sigma^{(2)}} \hookrightarrow \C^{\Sigma^+} \stackrel{\varrho}{\to} \C^\Gamma \to \C^{\Gamma''}
\]
is unital. If $\sigma \in \Sigma^{(2)}$ and $\gamma \in \Gamma''$, then $\varrho^{(3,2)}_\gamma(e_\sigma) \neq 0$ implies $k_\gamma \in \Jh(h_\sigma)$.
\alten

\altbn{\label{sd-1}}
$\Sd_J K$ is obtained inductively from the following Proposition. It says that, if we add the barycenter of some face of $|K|$ to the vertex set $V(K)$, we obtain a natural subdivision of $K$.

\begin{nprop}
Let $K$ be a simplicial complex and $\bar{\gamma} = \{\nu_0, \ldots,\nu_k\}$ a $k$-face of $K$. Then
\[
K_{\bar{\gamma}} := \{ \gamma \in K \, | \, \bar{\gamma} \not\subset \gamma \} \cup \{ \{\bar{\gamma}\} \cup (\gamma \setminus \{\nu_i\}) \, | \, \bar{\gamma} \subset \gamma \in K, \, i=0, \ldots,k\}
\]
defines a simplicial complex with vertex set $V(K_{\bar{\gamma}}) = V(K) \cup \{\bar{\gamma}\}$. Furthermore, the map $\beta: V(K_{\bar{\gamma}}) \to |K|$, where $\beta(\nu) := |\nu|$ for $\nu \in V(K)$ and $\beta(\bar{\gamma}) := \sum_{\nu \in \bar{\gamma}} \frac{1}{k+1} |\nu|$ is the barycenter of $\bar{\gamma}$ in $|K|$, induces a linear homeomorphism $\bar{\beta} : |K_{\bar{\gamma}}| \to |K|$. 
\end{nprop}
  
\begin{nproof}
It is easy to see that $K_{\bar{\gamma}}$ indeed defines a simplicial complex. If $\gamma =\{\nu_0, \ldots,\nu_l\}$ is a face of $K_{\bar{\gamma}}$, then (by definition of $K_{\bar{\gamma}}$ and $\beta$) $\{\beta(\nu_0),\ldots,\beta(\nu_l)\}$ is contained in some face of $|K|$, therefore $\beta$ induces a linear map $\bar{\beta} : |K_{\bar{\gamma}}| \to |K|$ by
\[
\bar{\beta}(t) = t_{\bar{\gamma}} \cdot \beta(\bar{\gamma}) + {\textstyle \sum_{\nu \in V(K)}} t_\nu \cdot |\nu|\, .
\]
If $s \in |K|$, set $d_s := \min \{s_\nu \, | \, \nu \in \bar{\gamma}\}$ and define $\alpha(s) \in |K_{\bar{\gamma}}|$ by 
\[
\alpha(s)_\nu := 
\left\{
\begin{array}{ll}
s_\nu & \mbox{ if } \nu \in V(K) \setminus \bar{\gamma}\\
s_\nu - d_s & \mbox{ if } \nu \in \bar{\gamma}\\
(k+1) d_s & \mbox{ if } \nu = \bar{\gamma}\, .
\end{array}
\right.
\]
Since the map $s \mapsto d_s$ is continuous, so is $\alpha : |K| \to |K_{\bar{\gamma}}|$; it is straightforward to check that $\alpha \circ \bar{\beta} = \id_{|K_{\bar{\gamma}}|}$ and that $\bar{\beta} \circ \alpha = \id_{|K|}$, hence $\bar{\beta}$ is a homeomorphism. 
\end{nproof}
\alten

\altbn{\label{sd-l}}
Let $\gamma_1, \ldots, \gamma_l$ be mutually distinct $k$-faces of $K$. Then $\gamma_2$ is a $k$-face of $K_{\gamma_1}$ ($\gamma_1 \neq \gamma_2$ are $k$-faces, so $\gamma_1 \not\subset \gamma_2$, hence $\gamma_2 \in K_{\gamma_1}$ by the definition of $K_{\gamma_1}$, cf.\ Proposition \ref{sd-1}). Therefore it makes sense to define $(K_{\gamma_1})_{\gamma_2}$ and, inductively, $K_m:= (\ldots(K_{\gamma_1})_{\gamma_2} \ldots )_{\gamma_m}$ for $m=1, \ldots,l$. Note that $|K_m| \approx |K|$ by the preceding Proposition. 

Now assume that $\gamma_i \cup \gamma_j \notin K$ whenever $i \neq j$. Under this condition, we check that $\{\gamma_i, \gamma_j\} \notin K_l $ for all $i \neq j$:

\noindent
Suppose that $\{\gamma_i,\gamma_j\} \in K_l$. We may assume $i < j$, then by the definition of the $K_m$, $\gamma_i$ is a vertex and $\gamma_j$ is a face of $K_{j-1}$; furthermore, we see that
\[
\{\gamma_i,\gamma_j\} \in K_l \Rightarrow \{\gamma_i,\gamma_j\} \in K_{l-1} \Rightarrow \ldots \Rightarrow \{\gamma_i,\gamma_j\} \in K_{j}\, .
\]
But then there must be some $\gamma \in K_{j-1}$ s.t.\ $\gamma_i$ (as a vertex of $K_{j-1}$) is an element of $\gamma$  and $\gamma_j$ (as a face of $K_{j-1}$) is a subset of $\gamma$; we may assume $\gamma = \gamma_j \cup \{\gamma_i\}$. Now by the definition of the $K_m$, $\gamma$ is a face of $K_m$ for $m= i, i+1, \ldots, j-1$. Therefore (this time by the definition of $K_i$) there is a face $\gamma'$ of $K_{i-1}$ containing both $\gamma_i$ and $\gamma_j$ as subsets. We may assume $\gamma' = \gamma_i \cup \gamma_j$. Since each vertex of $\gamma'$ is a vertex of $K$, again we see from the definition of the $K_m$, that $\gamma'$ in fact is a face of $K$, contradicting our assumption that $\gamma_i \cup \gamma_j \notin K$. Thus, $\{\gamma_i, \gamma_j\} \notin K_l$ if $i \neq j$.
\alten

\altbn
We are now prepared to construct the relative barycentric subdivision $\Sd_J K$. In the proof we will distinguish carefully between vertices $\nu \in V(K)$ of a simplicial complex and $0$-faces $\{\nu\} \in K$.

\begin{nproof} (of Proposition \ref{subdivision}) Suppose $h|_{\C^{\Sigma^{(1)}}}$ is $n$-decomposable w.r.t.\ the decomposition $\Sigma^{(1)} = \coprod_{j=0,\ldots,n} I_j$. \\
For $0 \le k \le n$ let $F^{(k)}$ consist of those $k$-faces $\gamma$ of $K$ for which $\gamma \subset (I_0 \cup \ldots \cup I_{k-1} \cup \Sigma^{(2)})$ (then, in particular, $* \notin \gamma$ for all $\gamma \in F^{(k)}$ and $F^{(0)} = \Sigma^{(2)}$).\\
If $\gamma = \{\nu_0, \ldots, \nu_k\} \in F^{(k)}$, then $h_{\nu_0} \cdots h_{\nu_k} \neq 0$ by the definition of $K$, so each $I_j$ contains at most one $\nu_i$. But this means that $\gamma \cap \Sigma^{(2)} \neq \emptyset$.\\  
If $\gamma_1, \ldots,\gamma_l$ are the distinct elements of $F^{(n)}$, define 
\[
K^{(n)} := (\ldots(K_{\gamma_1})_{\gamma_2} \ldots)_{\gamma_l}
\]
as in \ref{sd-l} and note that $|K^{(n)}| \approx |K|$ naturally. \\
Now let $\tau_1, \ldots, \tau_m$ be the distinct elements of $F^{(n-1)}$. None of the $\tau_i$ contains any of the $\gamma_j$, so we see from the definition of $K^{(n)}$ (cf.\ Proposition \ref{sd-1}), that the $\tau_i $ are faces of $K^{(n)}$. Set $K^{(n-1)} := (\ldots(K^{(n)}_{\tau_1})_{\tau_2} \ldots)_{\tau_m}$, then $|K^{(n-1)}| \approx |K^{(n)}| \approx |K|$. Proceed inductively to obtain $K^{(i)}$, $1\le i \le n$, s.t.\ $|K^{(1)}| \approx \ldots \approx |K^{(n)}| \approx |K|$. Define
\[
\Sd_J K := K^{(1)}\, ,
\]
then (i) and (ii) of \ref{subdivision} follow from our construction and Proposition \ref{sd-1}; note that
\[
\Gamma := V(\Sd_J K) = \Sigma^+ \cup F^{(1)} \cup \ldots \cup F^{(n)}\, .
\]
To obtain $K^{(1)}$ from $K$, we only changed faces of $K$ which intersect $\Sigma^{(2)}$ (namely the faces $F^{(0)} \cup \ldots \cup F^{(n)}$) and left the faces of $J$ invariant; this shows (iii). Furthermore, none of the faces in $F^{(0)} \cup \ldots \cup F^{(n)}$ contains the infinite point $*$, so if $t \in |K|$ is mapped to $|K^{(1)}|$, the $*$-coordinate of $\bar{\beta}^{-1}$ is equal to the $*$-coordinate of $t$. Therefore, $\bar{k}_{\{*\}} = \bar{h}_{\{*\}} \verk \bar{\beta}$, from which the third assertion of \ref{subdivision}(iv) follows.

Consider distinct elements $\gamma, \gamma'$ of $F^{(l)}$ and $\sigma \in I_l$ for some $0 \le l \le n$. Now if $\gamma \cup \gamma' \notin K$, then it is clear from our construction that $\gamma \cup \gamma' \notin K^{(k)}$ for any $k$ (otherwise, $\gamma \cup \gamma'$ and $F^{(m)}$ would have an element in common, but the elements of $\gamma \cup \gamma'$ are vertices, and not faces, of $K$). If $\gamma \cup \gamma' \in K$, then $\gamma \cup \gamma'$ is a $k$-face of $K$ for some $k >l$ and, since $\gamma \cup \gamma' \subset I_0 \cup \ldots \cup I_{l-1} \cup \Sigma^{(2)}$, we have $\gamma \cup \gamma' \in F^{(k)}$. Thus, again by construction, $\gamma \cup \gamma' $ is a vertex of $K^{(k)}$, but $\gamma \cup \gamma'$ (as a subset of $V(K)$) is not a face of $K^{(k)}$; more formally, we have $\gamma \cup \gamma' \in V(K^{(k)})$, hence $\{\gamma \cup \gamma'\} \in K^{(k)}$, but  $\gamma \cup \gamma' \notin K^{(k)}$.\\
Again it follows from the definition of the $K^{(m)}$, that $\gamma \cup \gamma'$ does not occur as a face of $K^{(m)}$ for any $m \le k$. \\
In particular, we have $\gamma \cup \gamma' \notin K^{(l)}$. Therefore, the elements of $F^{(l)}$ satisfy the condition of \ref{sd-l} (with $K^{(l+1)}$ in place of $K$), so \ref{sd-l} says that $\{\gamma,\gamma'\} \notin K^{(l)}$. By the same reasoning as above we see that $\{\gamma, \gamma'\}$ does not occur as a face of $K^{(m)}$ for any $m \le l$, thus $\{\gamma,\gamma'\} \notin \Sd_J K$.\\
Similarly, we show that $\{\gamma,\sigma\} \notin \Sd_J K$: Note that $\sigma \notin \gamma$ by definition of $F^{(l)}$. Now if $\{\gamma,\sigma\} \in \Sd_J K$, it follows from our construction  that $\gamma \cup \{\sigma\}$ must be a face of  $K^{(l+1)}$ and, inductively, that $\gamma \cup \{\sigma\}$ must be a face of $K$. On the other hand, we have $\gamma \in F^{(l)}$, $\sigma \in I_l$ and $\sigma \notin \gamma$ (so $\gamma \cup \{\sigma\}$ is an $(l+1)$-face of $K$ and $l<n$, since $K$ has at most $n$-faces), hence $\gamma \cup \{\sigma\} \in F^{(l+1)}$. But this in turn means that $\gamma \cup \{\sigma\}$ is a vertex, and not a face, of $K^{(l+1)}$, a contradiction, and we have $\{\gamma,\sigma\} \notin \Sd_J K$.

We are now prepared to show (iv) of Proposition \ref{subdivision}, namely that $\bar{k}|_{\C^{\Gamma \setminus \{*\}}}$ is $n$-decomposable.\\
For $0 \le j \le n$ define $I'_j := I_j \cup F^{(j)}$ (where $I_j$ comes from the decomposition of $\Sigma^{(1)}$), then $\Gamma \setminus \{*\} = \coprod_{j=0,\ldots,n} I'_j$.\\
If $\gamma \neq \gamma' \in I_j$, then $\bar{k}_\gamma \bar{k}_{\gamma'} = 0$: This is, because $\bar{\beta}$ maps the carriers of $\bar{k}_\gamma$ and $\bar{k}_{\gamma'}$ to the carriers of $\bar{h}_\gamma$ and $\bar{h}_{\gamma'}$, respectively, so $\bar{h}_\gamma \bar{h}_{\gamma'} = 0$ ($\bar{h}|_{\C^{I_j}}$ has strict order zero) implies $\bar{k}_\gamma \bar{k}_{\gamma'} = 0$.\\
If $\gamma \neq \gamma' \in F^{(j)}$ and $\sigma \in I_j$, then $\bar{k}_\gamma \bar{k}_{\gamma'} = \bar{k}_\gamma \bar{k}_{\sigma} = 0$, because $\{\gamma,\gamma'\}, \, \{\gamma,\sigma\} \notin \Sd_J K$, as we have seen before. \\
This shows that $\bar{k}|_{\C^{\Gamma \setminus \{*\}}}$ is $n$-decomposable w.r.t.\ the decomposition $\coprod_{j=0,\ldots,n} I'_j$. The second assertion of \ref{subdivision}(iv) follows directly from (iii), so the proof is complete.
\end{nproof}
\alten

\section{Proof of the main result}

\altbn
We are now prepared to prove the remaining part of Theorem \ref{rec.subh.}; in fact, we show a bit more.

\begin{ntheorem}
Let $A$ be a recursive subhomogeneous algebra of topological dimension not exceeding $n$. Then $A$ has a system $(F_\lambda, \psi_\lambda, \varphi_\lambda)_\Lambda$ of c.p.\ approximations with $\varphi_\lambda$ p.c.\ and $n$-decomposable $\forall \, \lambda$.
\end{ntheorem}

\begin{nproof}
If $A = M_r$, the theorem holds with the approximation $(M_r, \id_{M_r},\id_{M_r})$. So let $\Omega$ be compact and metrizable with $\dim \Omega \le n$ and let $r \in \N$. By induction we then have to show the following: Suppose $B$ is a recursive subhomogeneous algebra of topological dimension not exceeding $n$ for which the theorem holds, let $X \subset \Omega$ be closed and $\pi : B \to \Ch(X) \otimes M_r$ be a unital $*$-homomorphism. Then the assertion of the theorem holds for $A:= B \oplus_{\pi,X} (\Ch(\Omega) \otimes M_r)$.

{\it Step 1.} So let $\varepsilon > 0$ and $a_1, \ldots, a_k \in A_+$ with $\|a_l\| \le 1$ be given; we may assume $a_1 = \be_A$. We have to construct a c.p.\ approximation $(F,\psi,\varphi)$ s.t.\ $\|\varphi \psi(a_l) - a_l\| < \varepsilon$ and s.t.\ $\varphi$ is p.c.\ and $n$-decomposable. Set $b_l:= \beta(a_l) \in B$, $l=1, \ldots,k$, where $\beta: A \to B$ is the projection map.\\
Take $\alpha > 0$ s.t.\ $24(n+1) \alpha^\halb + 13 \alpha < \varepsilon$ and choose a c.p.c.\ approximation $(F'= \bigoplus_1^s M_{r_i}, \psi', \varphi')$ (of $B$) for $b_1, \ldots, b_k, b_1^2, \ldots, b_k^2$ within $\alpha$ such that $\varphi'$ is p.c.\ and $n$-decomposable. Then by Proposition \ref{pc-extension} there is a closed neighborhood $Y' \subset \Omega$ of $X$ and a c.p.c.\ map $\hat{\varphi} : F' \to B \oplus_{\pi, X} (\Ch(Y') \otimes M_r)$ s.t.\ $\hat{\varphi}$ is p.c., $n$-decomposable and
\begin{equation}{\label{(1)}}
\| \beta \verk \hat{\varphi}(x) - \varphi'(x) \| < \alpha \|x\| \; \forall \, x \in F'_+ .
\end{equation}
We may assume that for each $t \in Y'$ there is $\bar{t} \in X$ s.t.\ 
\begin{eqnarray}{\label{(2)}}
& &\| a_l(t) - a_l(\bar{t})\|, \, \| a_l^2(t) - a_l^2(\bar{t})\|, \, \| \hat{\varphi} \psi' (b_l) (t) - \hat{\varphi} \psi' (b_l) (\bar{t}) \| \, ,  \nonumber \\
& & \| \hat{\varphi} \psi' (b_l)^2 (t) - \hat{\varphi} \psi' (b_l)^2 (\bar{t}) \|\, , \, \| \hat{\varphi} \psi' (b_l^2) (t) - \hat{\varphi} \psi' (b_l^2) (\bar{t}) \|  <  \alpha \; \forall \, l \,.
\end{eqnarray}

{\it Step 2.} Since $\Omega$ is normal, there are open sets $V$, $W \subset \Omega$ and a closed set $Y \subset \Omega$ s.t.\ $X \subset W \subset Y \subset V \subset Y'$.\\
It is then straightforward to construct a finite collection $(U_\lambda)_\Lambda$ of open subsets of $\Omega \setminus X$ with the following properties:\\
(i) $(U_\lambda)_\Lambda$ is $n$-decomposable as a collection of subsets\\
(ii) $\Omega \setminus W \subset \bigcup_\Lambda U_\lambda$, $U_\lambda \cap (\Omega \setminus W) \neq \emptyset \; \forall \, \lambda$ and $U_\lambda \subset Y$ whenever $U_\lambda \cap W \neq \emptyset$\\
(iii) for each $\bar{\lambda} \in \Lambda$ there is $t_{\bar{\lambda}} \in U_{\bar{\lambda}}$ s.t.\ $t_{\bar{\lambda}} \notin \bigcup_{\Lambda \setminus \{\bar{\lambda}\}} U_\lambda$\\
(iv) $\| a_l(t) - a_l(t')\|$, $\|a_l^2(t) - a_l^2(t')\| < \alpha \, \forall \, t,t' \in U_\lambda$, $\lambda \in \Lambda$, $l=1, \ldots, k$\\
(v) $U_\lambda \subset Y' \; \forall \, \lambda \in \Lambda' := \{ \lambda \in \Lambda \, | \, U_\lambda \cap Y \neq \emptyset \}$\\
(vi) for each $\lambda \in \Lambda'$ there is $\mu(\lambda)$ with $\frac{1}{4(n+1)} < \mu(\lambda) < \frac{1}{2(n+1)}$ s.t., for all $i \in \{1, \ldots, s\}$, the projections
\[
q(\lambda, i):= \chi_{\mu(\lambda)} (\hat{\varphi}(\be_{M_{r_i}})|_{U_\lambda}) \in C^*(\hat{\varphi}(\be_{M_{r_i}})|_{U_\lambda}) \subset \Ch_b(U_\lambda) \otimes M_r
\]
exist (again, $\chi_\mu$ denotes the characteristic function of $[\mu,\infty)$).\\

Set $\Lambda^{(1)}:= \{ \lambda \in \Lambda \, | \, U_\lambda \cap W \neq \emptyset \} $ and $\Lambda^{(2)}:= \Lambda \setminus \Lambda^{(1)}$.

Next construct functions $g_\lambda \in \Ch_0(U_\lambda)$ for each $\lambda \in \Lambda$ with $ 0 \le g_\lambda \le 1$ s.t.\ $ g:= \sum_\Lambda g_\lambda$ satisfies $0 \le g \le 1$ and $g|_{\Omega\setminus W} \equiv 1$. We regard $g$ as an element of $A$ vanishing on $B$; note that $g \in \Zh(A)$ and that $(\be - g) g_\lambda =0$ for $\lambda \in \Lambda^{(2)}$.

{\it Step 3.} Now set $\bar{F} := \bar{F}^{(1)} \oplus \bar{F}^{(2)}$ with $\bar{F}^{(1)} := F'$, $\bar{F}^{(2)} := \bigoplus_{\Lambda^{(2)}} M_r$ and define $\bar{\psi} : A \to \bar{F}$ by $\psi' \verk \beta \oplus (\bigoplus_{\lambda \in \Lambda^{(2)}} \ev_{t_\lambda})$. The definition of $\bar{\varphi}: \bar{F} \to A$ requires some extra effort. For $\lambda \in \Lambda^{(2)}$ define $\bar{\varphi}^{(2)}_\lambda : M_r \to A$ by
\[
\bar{\varphi}^{(2)}_\lambda(x) := g_\lambda \cdot x \, ,
\]
where on the right hand side $x$ denotes the function $U_\lambda \to M_r$ with constant value $x$.\\
For the moment, let $\lambda \in \Lambda'$ be fixed. From (vi) we obtain projections $q(\lambda,i) \in \Ch_b(U_\lambda) \otimes M_r$ for $i=1, \ldots, s$. Note that
\[
\hat{\varphi}(\be_{M_{r_i}})|_{U_\lambda} - \frac{1}{2(n+1)} \le q(\lambda,i) \le 4(n+1) \cdot \hat{\varphi}(\be_{M_{r_i}})|_{U_\lambda} \, ,
\]
so 
\begin{eqnarray*}
\textstyle
\sum_{i=1}^s q(\lambda,i) & \ge & \textstyle \sum_{i=1}^s \hat{\varphi}(\be_{M_{r_i}})|_{U_\lambda} -  \halb \\
& \ge & \hat{\varphi} \psi' \beta (\be_A)|_{U_\lambda} - \halb \\
& \ge & \be|_{U_\lambda} - 2 \alpha - \halb \,,
\end{eqnarray*}
\noindent hence $\sum_{i=1}^s q(\lambda,i)$ is invertible in $\Ch_b(U_\lambda) \otimes M_r$. (We have used that $\hat{\varphi}$ is $n$-decomposable, so for any $t \in Y'$, $\hat{\varphi}(\be_{M_{r_i}})(t)$ is nonzero for at most $n+1$ indices $i$.) \\
We may thus apply Lemma \ref{l-piecewise-commuting} to obtain pairwise orthogonal projections $p(\lambda,i) \in C^*(q(\lambda,j) \, | \, j=1, \ldots , s) \subset \Ch_b(U_\lambda) \otimes M_r$, $i= 1, \ldots,s$, s.t.\
\begin{eqnarray}
{\textstyle \sum^s_{j=1} }p(\lambda,j) & = & \be|_{U_\lambda} \,, \nonumber\\
\left[ p(\lambda,i), \hat{\varphi}(F')|_{U_\lambda} \right] & = & 0 \,, {\label{(3)}}\\
p(\lambda,i) q(\lambda,i) & = & p(\lambda,i) \, . \nonumber 
\end{eqnarray} 
Note that (\ref{(3)}) in particular implies that
\begin{equation}{\label{(3a)}}
[p(\lambda,i)(t),q(\lambda',j)(t)] = [p(\lambda,i)(t),p(\lambda',j)(t)] = 0 
\end{equation}
for all  $\lambda,\lambda' \in \Lambda', \, t \in U_\lambda \cap U_{\lambda'}, \, i,j \in \{1,\ldots,s\}$.

We are now prepared to define $\bar{\varphi}^{(1)} : \bar{F}^{(1)} \to A$: \\
For $i \in \{ 1, \ldots,s\}$ and $x \in M_{r_i} \subset \bar{F}^{(1)}$ set
\begin{eqnarray*}
\textstyle
\lefteqn{\bar{\varphi}^{(1)}_i (x) := }\\
& &\textstyle (\be - g) \cdot \hat{\varphi}_i(x) + \sum_{\lambda \in \Lambda^{(1)}} g_\lambda \cdot p(\lambda,i) (q(\lambda,i) \hat{\varphi}(\be_{M_{r_i}})|_{U_\lambda})^{-1} \hat{\varphi}_i (x)|_{U_\lambda} \, ,
\end{eqnarray*}
where the inverses are taken in $q(\lambda,i) (\Ch_b(U_\lambda)\otimes M_r) q(\lambda,i)$.\\
Note that $\bar{\varphi}^{(1)}_i$ is well-defined and that
\begin{equation}{\label{(4)}}
\| (q(\lambda,i) \hat{\varphi}(\be_{M_{r_i}})|_{U_\lambda})^{-1} \| < 4(n+1) \; \forall \, i \in \{ 1,\ldots,s\}, \, \lambda \in \Lambda' \, .
\end{equation}
We shall need the following estimate later on: For all $\lambda \in \Lambda'$ and $t \in U_\lambda$ we have
\begin{eqnarray}
& &\| ({\textstyle \sum_i } \, p(\lambda,i) (q(\lambda,i) \hat{\varphi}(\be_{M_{r_i}})|_{U_\lambda})^{-1} \hat{\varphi}_i \psi_i' (b_l)) (t) - a_l(t) \|  \nonumber\\
&  & \; \; < 12 (n+1) \alpha^\halb + 4 \alpha \, .{\label{(5)}}
\end{eqnarray}
To see this, first note that there is $\bar{t} \in X$ s.t.\
\[
\begin{array}{rcccl}
b_l^2(\bar{t}) & \le & \varphi' \psi' (b_l)^2(\bar{t}) + 2 \alpha & \le & \hat{\varphi}\psi'(b_l)^2(\bar{t}) + 4 \alpha \\
& \stackrel{(\ref{(2)})}{\le} & \hat{\varphi} \psi' (b_l)^2(t) + 5 \alpha & \le & \hat{\varphi}(\psi'(b_l)^2)(t) + 5 \alpha \\
& \le & \hat{\varphi} \psi' (b_l^2) (t) + 5 \alpha & \stackrel{(\ref{(2)})}{\le} & \hat{\varphi}\psi'(b_l^2)(\bar{t}) + 6 \alpha \\
& \stackrel{(\ref{(1)})}{\le} & \varphi' \psi' (b_l^2)(\bar{t}) + 7 \alpha & \le & b_l^2 (\bar{t}) + 8 \alpha\,.
\end{array}
\]
As a consequence, 
\[
\| \hat{\varphi}(\psi'(b_l)^2)(t) - (\hat{\varphi} \psi'(b_l))^2 (t) \| \le 8 \alpha \, ;
\]
combining this with \cite{KW}, Lemma 3.5, we obtain
\[
\| \hat{\varphi}(x \psi'(b_l))(t) - (\hat{\varphi}(x) \hat{\varphi} \psi'(b_l)) (t) \| <  \| x \|  (8 \alpha)^\halb \; \forall \, x \in F' \, .
\]
Choosing $x = \be_{M_{r_i}}$ we get
\[
\| (p(\lambda,i) \hat{\varphi}_i \psi_i'(b_l))(t) - (p(\lambda,i) \hat{\varphi}(\be_{M_{r_i}}) \hat{\varphi} \psi'(b_l))(t) \| < 3 \alpha^\halb \, .
\]
But
\[
(p(\lambda,i) \hat{\varphi}(\be_{M_{r_i}})|_{U_\lambda})^{-1} = p(\lambda,i) (q(\lambda,i) \hat{\varphi}(\be_{M_{r_i}})|_{U_\lambda})^{-1} \le 4 (n+1) \cdot p(\lambda,i)
\]
by (\ref{(3)}) and (\ref{(4)}), so
\begin{eqnarray}
& & \| (p(\lambda,i) (q(\lambda,i) \hat{\varphi}(\be_{M_{r_i}})|_{U_\lambda})^{-1} \hat{\varphi}_i \psi'_i (b_l)) (t) - (p(\lambda,i) \hat{\varphi} \psi' (b_l))(t) \| \nonumber \\
& & \; \le 12 (n+1) \alpha^\halb \, ; \label{(5a)}
\end{eqnarray}
(\ref{(5)}) then follows by using (\ref{(1)}) and (\ref{(2)}) and the fact that the $p(\lambda,i)$ are pairwise orthogonal and sum up to $\be|_{U_\lambda}$.

The map $\bar{\varphi} := (\bar{\varphi}^{(1)}, \bar{\varphi}^{(2)}) : \bar{F} \to A$ is c.p.\ by construction; we check that it is contractive:
\begin{eqnarray*}
\textstyle
\bar{\varphi}(\be_{\bar{F}}) & = & \bar{\varphi}^{(1)}(\be_{\bar{F}^{(1)}}) + \bar{\varphi}^{(2)}(\be_{\bar{F}^{(2)}}) \\
& = & \textstyle (\be - g) \cdot \hat{\varphi}(\be_{F'}) + \sum^s_{i=1} \sum_{\lambda \in \Lambda^{(1)}} g_\lambda \cdot p(\lambda,i) + \sum_{\lambda \in \Lambda^{(2)}} g_\lambda \cdot \be_{M_r} \\
& \le & \textstyle (\be-g) \cdot \be_A + \sum_{\Lambda^{(1)}} g_\lambda \cdot \be_A + \sum_{\Lambda^{(2)}} g_\lambda \cdot \be_A \\
& = & \be_A \, .
\end{eqnarray*}
Furthermore, $\|\hat{\varphi}(\be_{F'})(t) - \be_A(t)\| \le 4 \alpha$ by (\ref{(1)}) and (\ref{(2)}), hence $\|\bar{\varphi}(\be_{\bar{F}}) - \be_A\| \le 4 \alpha$.\\
In fact $\bar{\varphi}$ has strict order zero on the summands of $F$ and is p.c., but for the moment it suffices to observe that $\bar{\varphi} \verk \bar{\iota} (\C^\Sigma)$ generates an abelian $C^*$-subalgebra of $A$, where $\Sigma := \{1, \ldots,s \} \stackrel{.}{\cup} \Lambda^{(2)}$ and $\bar{\iota} : \C^\Sigma \to \bar{F}$ is the canonical unital embedding. As it turns out, $\bar{\varphi} \verk \bar{\iota}|_{\C^\Sigma} $ satisfies a certain  order condition (cf.\ (\ref{(6)})). However, only the restriction of $\bar{\varphi} \verk \bar{\iota}$ to $\C^{\{1, \ldots,s\}}$ is $n$-decomposable, but we can use the idea of relative barycentric subdivision to obtain a modification $(F, \psi, \varphi)$ of $(\bar{F}, \bar{\psi}, \bar{\varphi})$ with the desired properties (in particular with $\varphi$ $n$-decomposable).

For $i \in \{1, \ldots,s\}$ we have
\[
\textstyle
\bar{\varphi}^{(1)}_i (\be_{M_{r_i}}) \le (\be - g) \cdot \hat{\varphi}^{(1)}_i (\be_{M_{r_i}}) + 4 (n+1) \sum_{\Lambda^{(1)}} g_\lambda \cdot \hat{\varphi}^{(1)}_i (\be_{M_{r_i}}) \, ,
\]
thus $\bar{\varphi}^{(1)}_i (\be_{M_{r_i}}) \perp \bar{\varphi}^{(1)}_j (\be_{M_{r_j}})$ if $\hat{\varphi}^{(1)}_i (\be_{M_{r_i}}) \perp \hat{\varphi}^{(1)}_j (\be_{M_{r_j}})$. But then the restriction $\bar{\varphi} \verk \bar{\iota}|_{\C^{\{1,\ldots,s\}}}$ is $n$-decomposable, since $\hat{\varphi}|_{\C^{\{1,\ldots,s\}}}$ is.

Now consider sets of mutually distinct elements $\{i_0, \ldots, i_m\} \subset \{1, \ldots,s\}$ and $\{\lambda_{m+1}, \ldots, \lambda_{n+1}\} \subset \Lambda^{(2)}$ for some $-1 \le m \le n$. By (ii) and the definition of $g$, $(\be-g) \cdot g_\lambda \equiv 0$ for $\lambda \in \Lambda^{(2)}$, hence $\bar{\varphi}^{(1)}_{i_0}(\be_{i_0}) \ldots \bar{\varphi}^{(2)}_{\lambda_{n+1}}(\be_{\lambda_{n+1}})$ is a sum of products of the form 
\[
g_{\lambda_0} \cdot p(\lambda_0,i_0) \ldots g_{\lambda_m} \cdot p(\lambda_m,i_m) g_{\lambda_{m+1}} \cdot \ldots \cdot g_{\lambda_{n+1}}
\]
with $\{\lambda_0, \ldots, \lambda_m\} \subset \Lambda^{(1)}$. Now if $0 \le l < l' \le m$ implies $\lambda_l \neq \lambda_{l'}$, then 
\[
g_{\lambda_0} \cdot \ldots \cdot g_{\lambda_m} \cdot g_{\lambda_{m+1}} \cdot \ldots \cdot g_{\lambda_{n+1}} = 0 \, ,
\]
because $(U_\lambda)_\Lambda$ is $n$-decomposable. If $\lambda_l = \lambda_{l'}$ for some $0 \le l < l' \le m$, then
\[
p(\lambda_l, i_l) p(\lambda_l, i_{l'}) =0 \, ,
\]
since for fixed $\lambda$ the $p(\lambda,i)$, $i=1, \ldots,s$, are pairwise orthogonal by construction. As a consequence, we obtain
\begin{equation}{\label{(6)}}
\bar{\varphi}_{i_0}^{(1)} (\be_{i_0}) \ldots \bar{\varphi}_{\lambda_{n+1}}^{(2)} (\be_{\lambda_{n+1}}) = 0 \, .
\end{equation}

{\it Step 4.} Now apply \ref{subdivision} with $C:= C^*(\bar{\varphi} \verk \bar{\iota} (\C^\Sigma), \be_A)$, $ \Sigma^{(1)}:= \{1, \ldots,s\}$ and $\Sigma^{(2)}:= \Lambda^{(2)}$; $h: \C^{\Sigma^+}  \to C$ is given by $h|_{\C^\Sigma}:= \bar{\varphi} \verk \bar{\iota}$ and $h|_{\C^{\{*\}}}:=\be_A - \bar{\varphi}(\be_{\bar{F}})$.  From \ref{properties-subdivision} we obtain an index set $\Gamma = \{*\} \stackrel{.}{\cup} \Sigma^{(1)} \stackrel{.}{\cup} \Gamma' \stackrel{.}{\cup} \Gamma''$ and a u.c.p.\ map $k: \C^\Gamma \to C$ s.t.\ $k|_{\C^{\Gamma \setminus \{*\}}}$ is $n$-decomposable and
\begin{eqnarray*}
k_i & \in & C^*(\bar{\varphi}_i(\be_{M_{r_i}}), \, g_\lambda \cdot\bar{\varphi_j}(\be_{M_{r_j}}),\,  g_\lambda, \, \be_A  \,|\, \lambda \in \Lambda^{(2)},\, j=1,\ldots,s ) \cap  \Jh(\bar{\varphi}_i(\be_{M_{r_i}})) \\
& & \mbox{for } i \in \{1, \ldots,s\} \, , \\
k_\gamma & \in & \Jh(g_{\nu(\gamma)}) \subset C \mbox{ for } \gamma \in \Gamma' \, , \\
k_\gamma & \in & C^*( g_\lambda, \, \be_A \,|\, \lambda \in \Lambda^{(2)} ) \cap \Jh(g_{\nu(\gamma)}) \mbox{ for } \gamma \in \Gamma'' \, ,
\end{eqnarray*}
where $\nu : \Gamma' \cap \Gamma'' \to \Lambda^{(2)}$ also comes from \ref{properties-subdivision} and $\Jh(M)$ denotes the ideal in $C$ generated by $M$. Using that $U_\lambda \cap W = \emptyset$ if $\lambda \in \Lambda^{(2)}$, it follows easily from \ref{subdivision}(iv) that  
\[
k_i(t) = \bar{\varphi_i}(\be_{M_{r_i}})(t) \; \forall \, t \in W\,.
\]
We have already seen that $\|\be_A - \bar{\varphi}(\be_{\bar{F}})\| \le 4 \alpha$.\\
For $\gamma \in \Gamma'$ we have $k_\gamma \cdot \sum_{i \in \Sigma^{(1)}} h_i \neq 0$ by the definition of $\Gamma'$. But then again $(\be - g) \cdot g_\lambda \equiv 0$ for $\lambda \in \Lambda^{(2)}$ implies that $g_{\nu(\gamma)} \cdot \sum_{\Lambda^{(1)}} g_\lambda \neq 0$. This in turn means that $U_{\nu(\gamma)} \cap W \neq \emptyset$, so $\nu(\gamma) \in \Lambda'$ (cf.\ (v)) and the projections $p(\nu(\gamma),i)$, $i=1, \ldots,s$, of Lemma \ref{l-piecewise-commuting} are well defined.\\
Note that $(\be -g) \cdot g_\lambda =0 \; \forall \, \lambda \in \Lambda^{(2)}$ also means that
\begin{eqnarray*}
\lefteqn{ \Jh(g_{\nu(\gamma)}) = }\\
& & C^*(g_\lambda, \, {\textstyle \sum_{\lambda' \in \Lambda^{(1)}}} \,  g_{\lambda'} \cdot p(\lambda',i) \,|\, i \in \{1, \ldots,s\}, \, \lambda \in \Lambda^{(2)}  ) \cap \Jh(g_{\nu(\gamma)}) \subset C \, .
\end{eqnarray*}
The following commutation relations are easily checked:
\begin{eqnarray*}
k_i & \in & (\bar{\varphi}_i(M_{r_i}))' \, , \, i=1, \ldots,s \, ,\\
k_\gamma & \in & (\hat{\varphi}_i(M_{r_i}) \cup \bar{\varphi}_i(M_{r_i}))' \, , \, \gamma \in \Gamma', \, i=1, \ldots,s \, ,\\
k_\gamma & \in & A' \, , \, \gamma \in \Gamma'' \, .
\end{eqnarray*}
From \ref{properties-subdivision} we also obtain a u.c.p.\ map $\varrho:\C^{\Sigma^+} \to \C^\Gamma$ s.t.\ $k \verk \varrho|_{\C^{\Gamma\setminus\{*\}}} = \bar{\varphi} \verk \bar{\iota}$, $\C^{ \{1, \ldots,s\} } \hookrightarrow \C^{\Sigma^+} \stackrel{\varrho}{\to} \C^\Gamma \to \C^{ \{1, \ldots,s\} }$ is the identity and $\varrho^{(3,2)} : \C^{\Sigma^{(2)}} \to \C^{\Gamma''}$ is unital.\\
If $\varrho^{(3,2)}_\gamma (e_\lambda) \neq 0$ for $\lambda \in \Lambda^{(2)}$ ($= \Sigma^{(2)}$) and $\gamma \in \Gamma''$, then we have $k_\gamma \in \Jh(g_\lambda)$, so $k_\gamma(t) \varrho^{(3,2)}_\gamma (e_\lambda) \neq 0$ implies $t \in U_\lambda$.

We finally turn to the definition of $(F, \psi, \varphi)$. Set 
\[
F := F^{(1)} \oplus F^{(2)} \oplus F^{(3)}
\]
with
\begin{eqnarray*}
F^{(1)} & := & {\textstyle \bigoplus_{i=1}^s } M_{r_i} = \bar{F}^{(1)} \, ,\\
F^{(2)} & := & {\textstyle \bigoplus_{\Gamma'}} ({\textstyle \bigoplus_{i=1}^s } M_{r_i}) \, ,\\
F^{(3)} & := & {\textstyle \bigoplus_{\Gamma''}} M_r \, .
\end{eqnarray*}
Write a map $\bar{\varrho} : \bar{F} \to F$ as a $3 \times 2$ matrix
\[
\bar{\varrho} := \left( \begin{array}{cc}
\id_{F^{(1)}} & 0 \\
{\textstyle \bigoplus_{\Gamma'}}\,  \id & 0 \\
0 & \varrho^{(3,2)} \otimes \id_{M_r}
\end{array}\right)
\]
and define
\[
\psi := \varrho \verk \bar{\psi} : A \to F \, ;
\]
$\psi$ clearly is c.p.c. Write $\varphi: F \to A$ as 
\[
\varphi := (\varphi^{(1)}, \, \varphi^{(2)}, \, \varphi^{(3)})
\]
with
\begin{eqnarray*}
\varphi^{(1)}_i (x) & := & k_i \cdot (\bar{\varphi}_i (\be_{M_{r_i}}))^{-1} \bar{\varphi}_i(x) \mbox{ for } i= 1, \ldots,s \mbox{ and } x \in M_{r_i} \, , \\
\varphi^{(2)}_{\gamma,i} (x) & := & k_\gamma \cdot p(\nu(\gamma),i) (q(\nu(\gamma),i) \hat{\varphi}_i(\be_{M_{r_i}})|_{U_{\nu(\gamma)}})^{-1} \hat{\varphi}_i(x) \mbox{ for } \gamma \in \Gamma', \\
& &   i \in \{ 1, \ldots, s \} \mbox{ and } x \in M_{r_i} \, , \\
\varphi^{(3)}_\gamma (x) & := & k_\gamma \cdot x \mbox{ for } \gamma \in \Gamma'' \mbox{ and } x \in M_r \, .
\end{eqnarray*}
Of course $\bar{\varphi}_i (\be_{M_{r_i}}) \in A$ need not be invertible; $(\bar{\varphi}_i (\be_{M_{r_i}}))^{-1} \bar{\varphi}_i (\, . \,)$ then stands for the $*$-homomorphism $: M_{r_i} \to A''$ associated to $\bar{\varphi}_i$ as in \cite{Wi1}, Lemma 1.1.3. Since $k_i \in \Jh (\bar{\varphi}_i (\be_{M_{r_i}}))$, a standard argument shows that $\varphi^{(1)}_i$ maps $M_{r_i}$ to $A$ and that $\varphi^{(1)}_i (\be_{M_{r_i}}) = k_i$.\\
The map $\varphi^{(2)}_{\gamma,i} : M_{r_i} \to A$ is well-defined, since $k_\gamma \in \Jh(g_{\nu(\gamma)}) \; \forall \, \gamma \in \Gamma'$.

{\it Step 5.} In the rest of the proof we check that $(F,\psi, \varphi)$ indeed is a c.p.\ approximation with the desired properties.

First note that $\varphi$ is contractive:
\begin{eqnarray*}
\varphi (\be_F) & = & {\textstyle \sum_{i=1}^s} \, k_i \cdot \be_A + {\textstyle \sum_{\gamma \in \Gamma'}} \, k_\gamma \cdot {\textstyle \sum_{i=1}^s} \, p(\nu(\gamma),i) + {\textstyle \sum_{\gamma \in \Gamma''}} \, k_\gamma \cdot \be_{M_r} \\
& \le & {\textstyle \sum_\Gamma} \, k_\gamma \cdot \be_A \\
& = & \be_A \, .
\end{eqnarray*}

$\varphi$ has strict order zero on the summands of $F$: Recall that
\begin{eqnarray*}
\lefteqn{\bar{\varphi}_i (x) := }\\
& &(\be - g) \cdot \hat{\varphi}_i(x) + {\textstyle \sum_{\lambda \in \Lambda^{(1)}}} \, g_\lambda \cdot p(\lambda,i) (q(\lambda,i) \hat{\varphi}(\be_{M_{r_i}})|_{U_\lambda})^{-1} \hat{\varphi}_i (x)|_{U_\lambda} \, ,
\end{eqnarray*}
and note that $(\be - g) \, , g_\lambda \, , p(\lambda,i) \, , q(\lambda,i)$ and $\hat{\varphi}_i(\be_{M_{r_i}})$ all commute with $\hat{\varphi}_i (M_{r_i})$; since $k_i$ commutes with $\bar{\varphi}_i(M_{r_i})$, $\ord \hat{\varphi}_i = 0$ implies that $\ord \varphi_i^{(1)}= 0 $ for $i \in \{1, \ldots,s\}$. A similar argument shows that $\ord \varphi_{\gamma,i}^{(2)} = 0$ for $\gamma \in \Gamma'$, $i \in \{1, \ldots,s\}$; obviously $\ord \varphi_\gamma^{(3)} = 0$ for $\gamma \in \Gamma''$.
 
$\varphi$ is p.c.: If $\gamma \in \Gamma''$, then $k_\gamma = \varphi^{(3)}_\gamma (\be_{M_r}) \in \Zh(A)$, so in particular 
\[
[\varphi^{(3)}_\gamma (\be_{M_r}), \, \varphi(F) ] = 0 \, .
\]
Next consider $i, \, j \in \{ 1, \ldots, s\}$; $\hat{\varphi}$ is p.c., hence w.l.o.g.\ we may assume that $[ \hat{\varphi}_i (\be_{M_{r_i}}), \, \hat{\varphi}_j (M_{r_j}) ] = 0$. For $\gamma \in \Gamma'$ we have
\begin{eqnarray*}
\lefteqn{ \varphi^{(1)}_i (\be_{M_{r_i}}) \varphi^{(2)}_{\gamma,j} (M_{r_j}) }\\
& = & k_i \cdot k_\gamma \cdot p(\nu(\gamma),j) (q(\nu(\gamma),j) \hat{\varphi}_j(\be_{M_{r_j}})|_{U_{\nu(\gamma)}})^{-1} \hat{\varphi}_j(M_{r_j}) \\
& = & k_\gamma \cdot p(\nu(\gamma),j) (q(\nu(\gamma),j) \hat{\varphi}_j(\be_{M_{r_j}})|_{U_{\nu(\gamma)}})^{-1} \hat{\varphi}_j(M_{r_j})
\cdot k_i \\
& = & \varphi^{(2)}_{\gamma,j} (M_{r_j}) \varphi^{(1)}_i (\be_{M_{r_i}}) \, ,
\end{eqnarray*}
because $k_i \cdot k_\gamma \in \Jh(g_{\nu(\gamma)}) \subset C^* (g_\lambda, \, {\textstyle \sum_{\lambda' \in \Lambda^{(1)}}} \, g_{\lambda'} \cdot p(\lambda',j), \,  | \, \lambda \in \Lambda^{(2)}, \, j =1,\ldots,s )$, which commutes with all other factors of the product (again we used that $(\be - g) \cdot g_{\nu(\gamma)} = 0$).\\
Similar arguments show that (for $\gamma, \, \gamma' \in \Gamma'$)
\begin{eqnarray*}
0 & = & [ \varphi^{(2)}_{\gamma,i} (\be_{M_{r_i}}), \, \varphi^{(1)}_j (M_{r_j}) ] \\
& = & [ \varphi^{(1)}_i (\be_{M_{r_i}}), \, \varphi^{(1)}_j (M_{r_j}) ] \\
& = & [ \varphi^{(2)}_{\gamma,i} (\be_{M_{r_i}}), \, \varphi^{(2)}_{\gamma',j} (M_{r_j}) ] \, . 
\end{eqnarray*}
It is then easy to define an order on $\{1, \ldots,s\} \stackrel{.}{\cup} (\stackrel{.}{\bigcup}_{\Gamma'} \{1, \ldots,s\}) \stackrel{.}{\cup} \Gamma''$ with respect to which $\varphi$ is piecewise commuting.

$\varphi$ is $n$-decomposable: Suppose $k|_{\C^{\Gamma \setminus \{*\}}}$ is $n$-decomposable w.r.t\ the decomposition $\Gamma \setminus \{*\} = \coprod_{j=0,\ldots,n} I_j$. Define a decomposition
\[
\textstyle
\{1,\ldots,s\} \cup (\Gamma' \times \{1,\ldots,s\}) \cup \Gamma'' = \coprod_{j=0,\ldots,n} I_j'
\]
by
\begin{eqnarray*}
I_j' \cap \{1,\ldots,s\} & := & I_j \cap \{1,\ldots,s\} \, , \\
I_j' \cap (\Gamma' \times \{1,\ldots,s\}) & := & (I_j \cap \Gamma') \times \{1,\ldots,s\} \, ,\\
I_j' \cap \Gamma'' & := & I_j \cap \Gamma'' \, .
\end{eqnarray*}
Then $\varphi$ is $n$-decomposable w.r.t.\ this decomposition $\coprod_{j=0,\ldots,n} I_j'$; this easily follows from the facts that $\varphi$ has strict order zero on the summands of $F$, that $k|_{\C^{\Gamma \setminus \{*\}}}$ is $n$-decomposable w.r.t.\ $\coprod_{j=0,\ldots,n} I_j$, and that $\varphi^{(2)}_{\gamma,i} (\be_{M_{r_i}}) \perp \varphi^{(2)}_{\gamma,j} (\be_{M_{r_j}})$ if $i \neq j$ (recall that $p(\lambda,i) \perp p(\lambda,j)$ for $\lambda \in \Lambda^{(2)}$ and $i \neq j$).

In the remainder of the proof we check that $(F, \psi, \varphi)$ indeed approximates the elements $a_1, \ldots, a_k$ within $\varepsilon$. It clearly suffices to do this separately for $B$ and for every $t \in \Omega \setminus X$; more precisely, we have to show that
\[
\| \beta ( a_l - \varphi \psi (a_l)) \|_B < \varepsilon
\]
and
\[
\| (a_l - \varphi \psi (a_l)) (t) \|_{M_r} < \varepsilon \; \forall \, t \in \Omega \setminus X \, .
\]
The former inequality is true, since
\begin{eqnarray*}
\| \beta(a_l) - \beta( \varphi \psi (a_l))\| & = & \| \beta(a_l) - \beta \hat{\varphi} \psi' \beta (a_l) \| \\
& \le & \| b_l - \varphi' \psi' (b_l)\| + \alpha \\
& < & 2 \alpha \\
& < & \varepsilon \, ;
\end{eqnarray*}
we check the latter. Note that
\[
\varphi \psi (a) = \varphi^{(1)} \psi^{(1)} (a) + \varphi^{(2)} \psi^{(2)} (a) + 
\varphi^{(3)} \psi^{(3)} (a) \, ,
\]
where $\psi^{(j)} : A \to F^{(j)}$ is $\psi$ followed by the projection onto $F^{(j)}$; we consider each summand $\varphi^{(j)} \psi^{(j)}$ separately.

We begin by showing that
\[
\| \varphi^{(1)} \psi^{(1)} (a_l)(t) - {\textstyle \sum_i} (k_i \cdot a_l)(t) \| < 12 (n+1) \alpha^\halb + 4 \alpha \, .
\]
If $t \notin Y$, then $\bar{\varphi}_i (\be_{M_{r_i}}) (t) = 0$ for $i= 1, \ldots, s$, so $k_i(t) = \varphi_i^{(1)}(\be_{M_{r_i}}) (t) = 0$ and there is nothing to show.\\ 
If $t \in Y \setminus W$, we obtain
\begin{eqnarray*}
\lefteqn{ \| \varphi^{(1)} \psi^{(1)} (a_l)(t) - {\textstyle \sum_i} k_i(t) \cdot a_l(t) \| } \\
& \stackrel{(\ref{(1)}),(\ref{(2)})}{\le} & \| \varphi^{(1)} \psi^{(1)} (a_l)(t) - {\textstyle \sum_i} k_i(t) \cdot \hat{\varphi} \psi' (b_l) (t) \| + 4\alpha \\
& = & \| {\textstyle \sum_i }\, k_i(t) \cdot (\bar{\varphi}_i (\be_{M_{r_i}}))^{-1} (t) \bar{\varphi}_i \psi' (b_l)(t) - {\textstyle \sum_i} \, k_i(t) \cdot \hat{\varphi} \psi' (b_l) (t) \| + 4\alpha \\
& = & \| {\textstyle \sum_i} \, k_i (t) \cdot ( {\textstyle \sum_{\Lambda^{(1)}}} \, g_\lambda \cdot p(\lambda,i))^{-1} (t)\\
& & \cdot ({\textstyle \sum_{\Lambda^{(1)}}} g_\lambda \cdot p(\lambda,i) (q(\lambda,i) \hat{\varphi}(\be_{M_{r_i}})|_{U_\lambda})^{-1} \hat{\varphi}_i \psi'_i (b_l))(t) \\
& & - {\textstyle \sum_i} \, k_i (t) \cdot ( {\textstyle \sum_{\Lambda^{(1)}}} \, g_\lambda \cdot p(\lambda,i))^{-1} (t) \cdot ({\textstyle \sum_{\Lambda^{(1)}}} g_\lambda \cdot p(\lambda,i) \hat{\varphi} \psi' (b_l))(t)\| \\
& & + 4\alpha \\  
& = & \| {\textstyle \sum_i} \, k_i (t) \cdot ( {\textstyle \sum_{\Lambda^{(1)}}} \, g_\lambda \cdot p(\lambda,i))^{-1} (t) \\
& & \cdot {\textstyle \sum_{\Lambda^{(1)}}} (g_\lambda \cdot p(\lambda,i)) (t) (p(\lambda,i) (q(\lambda,i) \hat{\varphi}(\be_{M_{r_i}})|_{U_\lambda})^{-1} \hat{\varphi}_i \psi'_i (b_l) \\
& & - p(\lambda,i) \hat{\varphi} \psi' (b_l)) (t) \| + 4\alpha \\
& \stackrel{(\ref{(5a)})}{\le} & \| {\textstyle \sum_i} \, k_i (t) \cdot ( {\textstyle \sum_{\Lambda^{(1)}}} \, g_\lambda \cdot p(\lambda,i))^{-1} (t) {\textstyle \sum_{\Lambda^{(1)}}} (g_\lambda \cdot p(\lambda,i)) (t)\| \cdot 12 (n+1) \alpha^\halb \\
& & + 4\alpha \\
& \le & 12( n+1) \alpha^\halb + 4\alpha \, .
\end{eqnarray*}
If $t \in W$, recall that $\bar{\varphi}_i (\be_{M_{r_i}}) (t) = k_i (t)$, so
\begin{eqnarray*}
\lefteqn{ \| \varphi^{(1)} \psi^{(1)} (a_l) (t) - {\textstyle \sum_i} k_i(t) \cdot a_l(t) \| } \\
& = & \| {\textstyle \sum_i} \, k_i(t) \cdot (\bar{\varphi}_i (\be_{M_{r_i}}))^{-1} (t) \bar{\varphi}_i \psi'_i (b_l) (t) - a_l (t) \| \\
& = & \| {\textstyle \sum_i} \, \bar{\varphi}_i \psi'_i (b_l) (t) - a_l (t) \| \\ 
& = & \| ( (\be - g) \cdot \hat{\varphi} \psi' (b_l) + {\textstyle \sum_i  \sum_{\Lambda^{(1)}}} \, g_\lambda \cdot p( \lambda,i) ( q(\lambda,i) \hat{\varphi} (\be_{M_{r_i}})|_{U_\lambda})^{-1} \\
& & - (\be -g) \cdot a_l + {\textstyle \sum_i \sum_{\Lambda^{(1)}}} \, g_\lambda \cdot p(\lambda,i) a_l )(t) \| \\
& \stackrel{(\ref{(5)}),(\ref{(5a)})}{\le} & (\be - g)(t) \cdot 4 \alpha + {\textstyle \sum_{\Lambda^{(1)}}} \, g_\lambda (t) \cdot (12(n+1) \alpha^\halb +4\alpha) \\
& < & 12 (n+1) \alpha^\halb + 4 \alpha \, .
\end{eqnarray*}
Next we check that 
\begin{eqnarray*}
\lefteqn{ \| \varphi^{(2)} \psi^{(2)} (a_l) - {\textstyle \sum_{\gamma \in \Gamma'}} \, k_\gamma \cdot a_l \| } \\
& = & \| {\textstyle \sum_{\gamma \in \Gamma'}} \, k_\gamma \cdot ( {\textstyle \sum_i} p(\nu(\gamma),i) (q(\nu(\gamma),i) \hat{\varphi} (\be_{M_{r_i}})|_{U_{\nu(\gamma)}})^{-1} \hat{\varphi}_i \verk \psi'_i(b_l) - a_l ) \| \\
& \le & (12 (n+1) \alpha^\halb + 4 \alpha) \cdot \| {\textstyle \sum_{\Gamma'}} \, k_\gamma \| \\
& \le & 12 (n+1) \alpha^\halb + 4 \alpha 
\end{eqnarray*}
and that
\begin{eqnarray*}
\lefteqn{ \| \varphi^{(3)} \psi^{(3)} (a_l) (t) - {\textstyle \sum_{\gamma \in \Gamma''}} \, (k_\gamma \cdot a_l) (t) \| } \\
& = & \| {\textstyle \sum_{\gamma \in \Gamma''}} \, k_\gamma (t) \cdot {\textstyle \sum_{\lambda \in \Lambda^{(2)}}} \, \rho_\gamma^{(3,2)} (e_\lambda) \cdot (a_l (t_\lambda) - a_l (t)) \| \\
& \le & {\textstyle \sum_{\Gamma'', \Lambda^{(2)}}} \, k_\gamma (t) \rho_\gamma^{(3,2)}(e_\lambda) \cdot \alpha \\
& \le & \alpha \, .
\end{eqnarray*}
Here we used that $\rho^{(3,2)}$ is unital (thus $\sum_{\Lambda^{(2)}} \, \rho_\gamma^{(3,2)} (e_\lambda) = 1$ for each $\gamma \in \Gamma''$) and that $t \in U_\lambda$, if $k_\gamma (t) \rho_\gamma^{(3,2)} (e_\lambda) \neq 0$.

As a consequence we obtain
\begin{eqnarray*}
\lefteqn{ \| \varphi \psi (a_l) - a_l \| } \\
& = & \| \varphi \psi (a_l) - {\textstyle \sum_\Gamma} k_\gamma a_l -(\be_A - \bar{\varphi}(\be_{\bar{F}})) a_l\| \\
& \le & \| \varphi^{(1)} \psi^{(1)} (a_l) - {\textstyle \sum_i} k_i \cdot a_l \| \\
& & + \| \varphi^{(2)} \psi^{(2)} (a_l) - {\textstyle \sum_{\Gamma'}} k_\gamma \cdot a_l \| \\
& & + \| \varphi^{(3)} \psi^{(3)} (a_l) - {\textstyle \sum_{\Gamma''}} k_\gamma \cdot a_l \| \\
& & + 4 \alpha\\
& < & 24 (n+1) \alpha^\halb + 13 \alpha \\
& < & \varepsilon \, , 
\end{eqnarray*}
so the proof is complete.
\end{nproof}
\alten





\end{document}